\documentclass[10pt,a4paper]{article}     
\usepackage{setspace,graphicx,amsmath,amssymb,graphics,float,fullpage,multirow,fancyhdr,amsthm,color,url,enumerate,subfig,pgf,tikz,mathrsfs,setspace,longtable}                       
\usepackage[font=small,labelfont=bf]{caption}
\usepackage[toc,page]{appendix}

\begin{document}             
\def\Pb{\mathbb{P}}
\def\Eb{\mathbb{E}}
\def\Rb{\mathbb{R}}
\def\Qb{\mathbb{Q}}
\def\Nb{\mathbb{N}}
\def\Zb{\mathbb{Z}}
\def\Tb{\mathbb{T}}
\def\Bb{\mathbb{B}}
\def\Tc{\mathcal{T}}
\def\Rc{\mathcal{R}}
\def\Sc{\mathcal{S}}
\def\Oc{\mathcal{O}}
\def\Uc{\mathcal{U}}
\def\Dc{\mathcal{D}}
\def\Jc{\mathcal{J}}
\def\Ec{\mathcal{E}}
\def\Cc{\mathcal{C}}
\def\Fc{\mathcal{F}}
\def\Lc{\mathcal{L}}
\def\Ac{\mathcal{A}}
\def\Gc{\mathcal{G}}
\def\Cc{\mathcal{C}}
\def\Mc{\mathcal{M}}
\def\Hc{\mathcal{H}}
\def\Bc{\mathcal{B}}
\def\Yc{\mathcal{Y}}
\def\Pc{\mathcal{P}}
\def\Ic{\mathcal{I}}
\def\Kc{\mathcal{K}}
\def\Pr{\mathbf{P}}
\def\Er{\mathbf{E}}
\def\Tk{\mathfrak{T}}
\def\Xk{\mathfrak{X}}
\def\Gk{\mathscr{G}}
\def\Ek{\mathscr{E}}
\def\Fs{\mathscr{F}}
\def\Ts{\mathscr{T}}
\def\Pt{\mathit{P}}
\def\Et{\mathit{E}}
\def\th{^{\text{th}}}
\def\oT{{\overline{\Tc}}}
\def\oR{{\overline{\rho}}}
\def\ooT{{\widetilde{\Tc}}}
\def\|{\\ \smallskip}
\def\d{\mathrm{d}}
\def\ind{\mathbf{1}}
\def\nin{n \rightarrow \infty}
\def\minf{m \rightarrow \infty}
\def\limn{\lim_{\nin}}
\def\limt{\lim_{t\rightarrow \infty}}
\def\limk{\lim_{k \rightarrow \infty}}
\def\limsn{\limsup_{\nin}}
\def\limin{\limsup_{\nin}}
\def\lime{\lim_{\varepsilon \rightarrow 0}}
\def\iid{\stackrel{\text{\tiny{i.i.d.}}}{\sim}}
\def\cd{\stackrel{\text{\tiny{d}}}{\rightarrow}}
\def\cl{\stackrel{\Lc}{\rightarrow}}
\def\ed{\stackrel{\text{\tiny{d}}}{=}}
\def\qqquad{\qquad \quad}
\def\qqqquad{\qqquad \quad}
\def\qqqqquad{\qqqquad \quad}
\def\qqqqqquad{\qqqqquad \quad}
\def\qqqqqqquad{\qqqqqquad \quad}
\usetikzlibrary{arrows}
\newtheorem{dfn}{Definition}
\newtheorem{thm}{Theorem}
\newtheorem{lem}{Lemma}[section]
\newtheorem{prp}[lem]{Proposition}
\newtheorem{cly}[lem]{Corollary}
\newtheorem{conj}[lem]{Conjecture}
\newtheorem{fig}{Figure}
\pagenumbering{arabic} 
\singlespacing
\allowdisplaybreaks
\title{Central limit theorems for biased randomly trapped random walks on $\Zb$}
\author{Adam Bowditch, University of Warwick}
\date{}
\maketitle
\numberwithin{equation}{section}
\begin{abstract}
We prove CLTs for biased randomly trapped random walks in one dimension. 
In particular, we will establish an annealed invariance principal by considering a sequence of regeneration times under the assumption that the trapping times have finite second moment.
In a quenched environment, an environment dependent centring is determined which is necessary to achieve a central limit theorem. 
As our main motivation, we apply these results to biased walks on subcritical Galton-Watson trees conditioned to survive and prove a tight bound on the bias required to obtain such limiting behaviour. 
\end{abstract}

\let\thefootnote\relax\footnote{\textit{MSC2010 subject classifications:} Primary 60K37, 60F05, 60F17; secondary 60J80. }

\section{Introduction}\label{Int}
In this paper, we investigate biased randomly trapped random walks (RTRWs) on $\Zb$ and apply the results to subcritical Galton-Watson trees conditioned to survive. Randomly trapped random walks were first introduced in \cite{arcacero} where it is shown that the possible scaling limits belong to a certain class of time changed Brownian motions. The purpose of the RTRW is to generalise models such as the Bouchaud trap model (see \cite{bouc}, \cite{foisne} and \cite{zi}) and provide a framework for studying random walks on other random graphs in which trapping naturally occurs, such as biased random walks on percolation clusters (see  \cite{dh}, \cite{frha} and \cite{szit}) and random walk in random environment (see \cite{kekosp} and \cite{szze}). Higher dimensional ($d\geq 2$) unbiased randomly trapped random walks have been studied further in \cite{cewa} where a complete classification of the possible scaling limits is given. In recent years there has been much progress in models involving trapping phenomena; a review of recent developments in a range of models of directionally transient and reversible random walks on underlying graphs such as supercritical GW-trees and supercritical percolation clusters is given in \cite{arfr}.

A subcritical GW-tree conditioned to survive consists of a semi-infinite path (called the backbone) with GW-trees as leaves. Typically the leaves are quite short therefore the walk on the tree does not deviate too far from the backbone. For this reason we have that the walk on the tree behaves like a randomly trapped random walk on $\Zb$ with holding times distributed as excursion times in GW-trees. Biased walks on subcritical GW-trees are, therefore, a natural example of the randomly trapped random walk. Furthermore, they exhibit interesting behaviour as the relationship between the bias and offspring law influences the trapping. In this paper we are only concerned with ballistic walks. We note that the sub-ballistic regimes for the biased walk on the subcritical tree have been studied in \cite{bowd} where it is shown that either a strong bias or heavy tails of the offspring law can slow the walk into a sub-ballistic phase. 

Critical and supercritical GW-trees have also received much attention. In \cite{crfrku}, it is shown that the walk on the critical GW-tree conditioned to survive is always sub-ballistic; this is studied further and it is shown that the walk belongs to the universality class of one-dimensional trapping models with slowly-varying tails. For the biased walk on the supercritical GW-tree, it is shown in \cite{lypepe} that when the bias is small the walk is recurrent, when the bias is large enough the walk is sub-ballistic and there is some intermediate range for the bias such that the walk is ballistic. The ballistic phase for this walk is studied further in \cite{ai} where an expression of the speed is given and in \cite{arfrgaha} where appropriate scaling sequences for the sub-ballistic phase are shown. The traps formed in the supercritical tree resemble those in the subcritical tree and it has been shown in \cite{bowd} that the walks observe similar scaling regimes. 

A technique is developed in \cite{bosz} which can be used to extend an annealed invariance principal to a quenched result. This is applied in \cite{peze} to prove a quenched functional central limit theorem for the walk on the supercritical tree when the offspring distribution has exponential moments and no deaths. The condition of exponential moments is purely technical. However, because the offspring law has no deaths, the supercritical tree does not have traps which represents a significant simplification of the problem. Due to the similarity of the traps in the supercritical and subcritical GW-trees with leaves, a key motivation of this paper is to be able to extend the result of \cite{peze} to allow for deaths in the offspring law. This will be done in \cite{bopre}.

The argument of \cite{bosz} relies on multiple copies of the walk eventually having separate escape paths; this creates a mixing of the environment which removes the dependency on the specific environment. In the one-dimensional model, the walk is forced to escape along a single route and therefore visits every trap in the positive segment of the environment. This results in the walk accumulating environment dependent fluctuations. In Section \ref{QuCLT} we prove a quenched central limit theorem with environment dependent centring which demonstrates why this occurs in further detail.

We next introduce the models of interest and results in more detail. We consider the randomly trapped random walk model in which the underlying walk $(Y_k)_{k\geq 0}$ is a simple, biased random walk on $\Zb$. That is, we write $Y_k:=\sum_{j=1}^k\chi_j$ for a sequence of i.i.d.\ random variables $(\chi_j)_{j\geq 1}$  satisfying $\Pt(\chi_j=-1)=(\beta+1)^{-1}=1-\Pt(\chi_j=1)$. For $x \in \Zb$ write 
\[L(x,n):=\sum_{k=0}^n\ind_{\{Y_k=x\}}\]
for the local time of $Y$ at site $x$ by time $n$. A random environment $\omega$ is a sequence of $(0,\infty)$-valued probability measures $(\omega_x)_{x \in \Zb}$ with environment law $\Pr:=\pi^{\otimes \Zb}$ for some fixed law $\pi$. For a fixed environment $\omega$, let $(\eta_{x,i})_{x \in \Zb, i\geq 0}$ be independent with $\eta_{x,i} \sim \omega_x$. Writing 
\begin{flalign*}
 S_n \; := \; \sum_{x \in \Zb} \sum_{i=1}^{L(x,n-1)}\eta_{x,i} \; = \; \sum_{k=0}^{n-1}\eta_{Y_k,L(Y_k,k)} \qquad \text{ and } \qquad S^{-1}_t=\inf\{k\geq 0:S_k>t\}
\end{flalign*}
we then define the randomly trapped random walk by
\[X_t:=Y_{S^{-1}_t}.\]
This process is then a continuous time random walk on $\Zb$ with $k^{th}$ holding time $\eta_k:=\eta_{Y_k,L(Y_k,k)}$ and we write $\eta:=(\eta_k)_{k\geq 0}$ to be the sequence of holding times. For convenience we will define $S_t=S_{\lfloor t\rfloor}$ where $\lfloor t\rfloor:=\max\{k \in \Zb:k\leq t\}$ for non-integer $t \in \Rb$.  Let $\Pt^\omega$ to denote the law over $X$ for fixed environment $\omega$ and $\Pb(\cdot)=\int \Pt^\omega(\cdot) \Pr(\d \omega)$ the annealed law. 

In Section \ref{AnCLT} we prove a functional law of large numbers and central limit theorem for the randomly trapped random walk on $\Zb$. That is, in Corollary \ref{LinX}, we show that if the bias is positive ($\beta>1$) and the expected holding time $\eta_0$ (under $\Pb$) is finite, then $X_{nt}/n$ converges to the deterministic process $\nu t$ where $\nu$ is a known constant (called the speed). In Corollary \ref{EIN} we then show that this speed satisfies an Einstein relation; more specifically, the derivative of the speed with respect to the bias approaches half of the diffusion coefficient of the unbiased walk as $\beta \rightarrow 1^+$. The main result of the section is Theorem \ref{AIP} which states that, if $\beta>1$ and $\Eb[\eta_0^2]<\infty$ then for some $\varsigma \in (0,\infty)$
\[B_t^n:=\frac{X_{nt}-nt\nu}{\varsigma \sqrt{n}}\]
converges in distribution to a Brownian motion. We prove this by considering a renewal argument similar to that of \cite{sz}. 

In Section \ref{QuCLT} we adapt the technique used in \cite{go} (to prove a quenched CLT for a random walk in random environment) to derive a quenched central limit theorem with an environment dependent centring for the RTRW (Theorem \ref{QCLT}). That is, for a fixed environment $\omega$, we show that if $\beta>1$, $\Eb[\eta_0^2]<\infty$ and for some $\varepsilon>0$ that $\Er[\Et^\omega[\eta_0]^{2+\varepsilon}]<\infty$ then there exists an environment dependent centring $\Gc^\omega(n)$ such that for some $\vartheta \in (0,\infty)$
\[\frac{X_{n}-\Gc^\omega(n)}{\vartheta \sqrt{n}}\]
converges in distribution to a standard Gaussian. In particular, we show that the known function $\Gc^\omega(n)$ can be written as the annealed, deterministic centring with an environment dependent correction where this correction is a sum of centred i.i.d.\ random variables with non-zero variance under the environment law. This shows that the correction obeys a central limit theorem under $\Pr$, thus has $\sqrt{n}$ fluctuations and is, therefore, necessary.

In Section \ref{Sub} we then apply these results to the biased random walk on the subcritical GW-tree conditioned to survive. Let $\xi$ denote the offspring distribution of a GW-process with mean $\mu \in (0,1)$ and variance $\sigma^2<\infty$ then let $|X_n|$ denote the graph distance between the walk at time $n$ and the root of the tree. In Lemma \ref{spd} and Corollary \ref{SubTreeSpd} we show that if $\beta \in (1,\mu^{-1})$ then the expected trapping time is finite and the law of large numbers holds for $|X_n|$ with speed
\[\nu_\beta=\frac{\mu(\beta-1)(1-\beta\mu)}{\mu(\beta+1)(1-\beta\mu)+2\beta(\sigma^2-\mu(1-\mu))}.\]
By Lemma \ref{spd} the expected trapping time is infinite when $\sigma^2=\infty$ or $\beta\geq \mu^{-1}$. Moreover, when $\beta<1$ the walk is recurrent. It therefore follows that the conditions $\sigma^2<\infty$ and $\beta \in (1,\mu^{-1})$ are also necessary for the walk the be ballistic.

We then develop the speed result to an Einstein relation in Lemma \ref{SubEIN}. In Corollary \ref{AnnSecMom} we prove that if $\Er[\xi^3]<\infty$ and $\beta \in (1,\mu^{-1/2})$, then
\[B_t^n:=\frac{|X_{nt}|-nt\nu_\beta}{\varsigma \sqrt{n}}\]
converges in distribution, under the annealed law $\Pb$, to a Brownian motion. In Proposition \ref{enqusq} we prove the quenched analogue. That is, if we also have that $\Er[\xi^{3+\varepsilon}]<\infty$ for some $\varepsilon>0$ then for $\Pr$-a.e.\ tree $\Tc$,
\[\frac{|X_{n}|-\Gc^\Tc(n)}{\vartheta \sqrt{n}}\]
converges in distribution to a standard Gaussian under $\Pt^\Tc$ where $\Gc^\Tc$ is an environment dependent centring. Moreover, under these conditions, $(\nu_\beta n-\Gc^\Tc(n))/\sqrt{n}$ converges in distribution to a Gaussian random variable (with strictly positive variance) under $\Pr$, which confirms the need for an environment dependent centring in the quenched CLT for the biased random walk on the subcritical GW-tree conditioned to survive. We also show, in Lemma \ref{Nec}, that the conditions $\beta^2\mu< 1$ and $\Er[\xi^3]<\infty$ are necessary for the previous results.

Let $\beta_c$ denote the smallest bias such that the walk is sub-ballistic whenever $\beta>\beta_c$; then, note that $\beta_c=\mu^{-1}$ for the walk on the subcritical GW-tree. The necessity that $\beta<\mu^{-1/2}$ from Lemma \ref{Nec} supports \cite[Conjecture 3.1]{arfr} which states that a quenched central limit should hold on the supercritical tree only when $\beta<\beta_c^{1/2}$ which, as noted above, we will prove in an accompanying paper.

\section{A functional law of large numbers and central limit theorem}\label{AnCLT}
The main aim of this section is to derive an annealed functional central limit theorem for the RTRW model with positive bias. That is, we show that $(X_{nt}-nt\nu(\beta))n^{-\frac{1}{2}}$ converges in law to a Brownian motion under a second moment condition by considering a regeneration argument similar to that used in \cite{sz}.

We first prove a straightforward speed result for the walk. By using ergodicity of the sequence of holding times we show that the sum of the first $nt$ holding times converges almost surely, and then extend this to a result for the walk. Notice that we consider the unbiased case ($\beta=1$) as well as the positive bias case. This will be used when we prove an Einstein relation for the walk. The key to Proposition \ref{Sconv} is \cite[Lemma 2.1]{cewa} which states that the left shift on sequences ($\theta(\eta_0,\eta_1,...)=(\eta_1,\eta_2,...)$) acts ergodically on $\eta$ under $\Pb$. This holds for any non-degenerate random walk on a fixed environment with i.i.d.\ traps which is why we we can extend the following result to the unbiased case.
\begin{prp}\label{Sconv}
 Suppose $\beta\geq 1$ and that $\Eb[\eta_0]<\infty$ then $S_{nt}/n$ and $S^{-1}_{nt}/n$ converge $\Pb$-a.s.\ on $D([0,\infty),\Rb)$ endowed with the Skorohod metric to the deterministic processes $\Sc_t=\Eb[\eta_0]t$ and $\Sc_t^{-1}=\Eb[\eta_0]^{-1}t$ respectively.
 \begin{proof}
 We begin by showing convergence of $S_{nt}/n$. Since $\theta$ acts ergodically on $\eta$ under $\Pb$ and $\Eb[\eta_0]<\infty$ we have that $f(\eta):=\eta_0$ is integrable. Therefore, by the ergodic theorem
 \begin{flalign*}
  \limn \frac{S_{nt}}{n} \; = \; \limn \frac{1}{n} \sum_{k=0}^{\lfloor nt\rfloor -1} \eta_k \; = \; \limn t \frac{1}{nt}\sum_{k=0}^{\lfloor nt\rfloor -1} f(\theta^k\eta) \; = \; t\Eb[f(\eta)]
 \end{flalign*}
almost surely. The sequence of functions $S_{nt}/n$ are increasing in $t$ and the limit $\Sc_t=\Eb[\eta_0]t$ is continuous therefore the convergence holds uniformly over $t \in [0,T]$ for $T<\infty$.

Since $\Sc_t$ is strictly increasing we have the desired convergence of $S^{-1}_{nt}/n$ by \cite[Corollary 13.6.4]{wh}.

%
 \end{proof}
\end{prp}

Using Proposition \ref{Sconv}, we are now able to show the speed result for the randomly trapped random walk.
\begin{cly}\label{LinX}
 Suppose $\beta>1$ and that $\Eb[\eta_0]<\infty$, then $X_{nt}/n$ converges $\Pb$-a.s.\ on $D([0,\infty),\Rb)$ endowed with the Skorohod metric to the process $\nu t$ where \[ \nu := \frac{(\beta-1)}{\Eb[\eta_0](\beta+1)}\]
 uniformly over $t$ in compact sets.
 \begin{proof}
 Notice that 
 \[\frac{X_{nt}}{n} \; = \; \frac{Y_{S^{-1}_{\lfloor nt\rfloor}}}{n} \; = \;  \frac{Y_{S^{-1}_{\lfloor nt\rfloor}}}{S^{-1}_{\lfloor nt\rfloor}}\cdot\frac{S^{-1}_{\lfloor nt\rfloor}}{n}. \]
  By the law of large numbers $n^{-1}Y_n$ converges a.s.\ to $(\beta-1)/(\beta+1)$ therefore by Proposition \ref{Sconv} and using that $\Sc_t^{-1}$ is continuous we indeed have the desired result.
 \end{proof}

\end{cly}

An additional result that can be deduced from Proposition \ref{Sconv} and Corollary \ref{LinX} is that the following Einstein relation holds. 
\begin{cly}\label{EIN}
 Suppose $\Eb[\eta_0]<\infty$. The unbiased $(\beta=1)$ walk $X_{\lfloor nt \rfloor}n^{-1/2}$ converges in $\Pb$-distribution on $D([0,\infty),\Rb)$ endowed with the Skorohod metric to a scaled Brownian motion with variance $\Upsilon=\Eb[\eta_0]^{-1}$. Moreover, 
\[\lim_{\beta\rightarrow 1^+}\frac{\nu}{\beta-1}=\frac{\Upsilon}{2}\]
where $\nu$ is the speed calculated in Corollary \ref{LinX} for the $\beta$-biased walk.
\begin{proof}
For $\beta=1$ we have that $Y_{nt}$ is the sum of i.i.d.\ copies of the random variable $\chi$ satisfying $\Pt(\chi=1)=1/2=\Pt(\chi=-1)$ thus by Donsker's invariance principle $Y_{nt}n^{-1/2}$ converges in $\Pb$-distribution on $D([0,\infty),\Rb)$ endowed with the Skorohod metric to a standard Brownian motion. 

By Proposition \ref{Sconv} we have that $S^{-1}_{\lfloor nt\rfloor}/n$ converges $\Pb$-a.s.\ to the deterministic process $t/\Eb[\eta_0]$ uniformly over $t\leq T$. By continuity of the limiting Brownian motion and \cite[Theorem 13.2.1]{wh}, we have that
 \begin{flalign*}
  \frac{X_{nt}}{\sqrt{n}}  = \frac{Y_{S^{-1}_{\lfloor nt\rfloor}}}{\sqrt{n}} \qquad \text{ and } \qquad  \frac{Y_{nt/\Eb[\eta_0]}}{\sqrt{n}}
 \end{flalign*}
converge to the same limiting process, which is a scaled Brownian motion with variance $\Upsilon=\Eb[\eta_0]^{-1}$.

Moreover, by Corollary \ref{LinX} we have that, for $\beta>0$, \[ \nu = \frac{(\beta-1)}{\Eb[\eta_0](\beta+1)},\]
and therefore we indeed have that \[\lim_{\beta\rightarrow 1^+}\frac{\nu}{\beta-1}=\frac{\Upsilon}{2}.\]
\end{proof}
\end{cly}

We now move on to proving an annealed functional CLT which is the main result of this section. That is, we show that
\[B_t^n:=\frac{X_{nt} -nt\nu}{\varsigma \sqrt{n}}\] 
converges in $\Pb$-distribution to a standard Brownian motion for some $\varsigma^2 \in (0,\infty)$. We want to write $X_{nt}$ as approximately the sum of i.i.d.\ centred random variables with finite second moments. Let $\kappa_0=0$ and, for $j=1,2,...$, define $\kappa_j:=\inf\{m>\kappa_{j-1}: \; \{Y_l\}_{l=0}^{m-1} \cap \{Y_l\}_{l=m}^\infty = \phi \}$ to be the regeneration times of the walk $Y$. We then have that $S_{\kappa_j}$ for $j\geq 1$ are regeneration times for $X$ and we write 
\[Z_j\;:=\;\left(X_{S_{\kappa_j}}-X_{S_{\kappa_{j-1}}}-\left(S_{\kappa_j}-S_{\kappa_{j-1}}\right)\nu\right)\;=\;\left(Y_{\kappa_j}-Y_{\kappa_{j-1}}-\left(S_{\kappa_j}-S_{\kappa_{j-1}}\right)\nu\right).\]

An important result of \cite[Lemma 5.1]{depeze} and the remark leading to it is that the time between regenerations of $Y$ has exponential moments, that is
\begin{flalign}\label{expMom}
\Pb(\kappa_{j+1}-\kappa_j>n) \leq Ce^{-cn}
\end{flalign}
for any $j\geq 1$ and some constants $C,c$. Moreover, $Y_{\kappa_{j+1}}-Y_{\kappa_j}\leq \kappa_{j+1}-\kappa_j$ therefore the distance between regeneration points also has exponential moments. 

\begin{lem}\label{Zcent}
Suppose that $\beta>1$ and $\Eb[\eta_0]<\infty$ then $\{Z_j\}_{j\geq 2}$ are centred and i.i.d.\ under $\Pb$.
\begin{proof}
By \cite{depeze} we have that the sections of the walk $(Y_{i+\kappa_j}-Y_{\kappa_j})_{i=1}^{\kappa_{j+1}-\kappa_j}$, for $j\geq 1$, are i.i.d.\ therefore, since the traps $(\omega_x)_{x \in\Zb}$ are i.i.d.\ and independent of the walk $Y$, we have that the sequences $(\eta_k)_{k=\kappa_j}^{\kappa_{j+1}-1}$ are i.i.d. It follows that $\{Z_j\}_{j\geq 2}$ are i.i.d.\ under $\Pb$. 

It remains to show that $Z_j$ are centred. Since the distribution of a given holding time is independent of the regeneration times of $Y$ and $\Eb[\eta_0]<\infty$ we have that
\begin{flalign}\label{nuS}
\nu\Eb[S_{\kappa_2}-S_{\kappa_{1}}] = \frac{\beta-1}{(\beta+1)\Eb[\eta_0]}\Eb\left[\sum_{k=\kappa_1}^{\kappa_2-1}\Eb[\eta_k|\kappa_2,\kappa_1]\right]= \frac{\beta-1}{\beta+1}\Eb[\kappa_2-\kappa_1].
\end{flalign}
We want to show this is equal to $\Eb[X_{S_{\kappa_j}}-X_{S_{\kappa_{j-1}}}]=\Eb[Y_{\kappa_2}-Y_{\kappa_1}]$. By (\ref{expMom}) the time between regenerations and distance between regeneration points have exponential moments hence, by the law of large numbers,
\begin{flalign}\label{YkappaConFrac}
 &&\frac{\sum_{j=2}^mY_{\kappa_j}-Y_{\kappa_{j-1}}}{m} & \rightarrow \Eb[Y_{\kappa_2}-Y_{\kappa_1}],& \notag\\
 &&\frac{\sum_{j=2}^m\kappa_j-\kappa_{j-1}}{m} & \rightarrow \Eb[\kappa_2-\kappa_1] &\notag\\
\text{and therefore} &&& \notag\\
&&\frac{\sum_{j=2}^mY_{\kappa_j}-Y_{\kappa_{j-1}}}{\sum_{j=2}^m\kappa_j-\kappa_{j-1}} & \rightarrow \frac{\Eb[Y_{\kappa_2}-Y_{\kappa_1}]}{\Eb[\kappa_2-\kappa_1]} &
\end{flalign}
$\Pb$-a.s. as $\minf$. However, 
\[\frac{\sum_{j=2}^mY_{\kappa_j}-Y_{\kappa_{j-1}}}{\sum_{j=2}^m\kappa_j-\kappa_{j-1}} = \frac{Y_{\kappa_m}}{\kappa_m}\left(1+\frac{\kappa_1}{\kappa_m-\kappa_1}\right)-\frac{Y_{\kappa_1}}{\kappa_m-\kappa_1}\]
where $Y_{\kappa_1}/(\kappa_m-\kappa_1), \; \kappa_1/(\kappa_m-\kappa_1)$ converge $\Pb$-a.s.\ to $0$. Furthermore, by the law of large numbers, $Y_{\kappa_m}/\kappa_m$ converges $\Pb$-a.s.\ to $(\beta-1)/(\beta+1)$ therefore
\begin{flalign}\label{YKappa}
\frac{\sum_{j=2}^mY_{\kappa_j}-Y_{\kappa_{j-1}}}{\sum_{j=2}^m\kappa_j-\kappa_{j-1}}\rightarrow \frac{\beta-1}{\beta+1}.
\end{flalign}
By (\ref{nuS}), (\ref{YkappaConFrac}) and (\ref{YKappa}) we then have that
\[\Eb[Y_{\kappa_2}-Y_{\kappa_1}]=\frac{\beta-1}{\beta+1}\Eb[\kappa_2-\kappa_1]=\nu\Eb[S_{\kappa_2}-S_{\kappa_{1}}].\]
Therefore $Z_j$ are centred, as desired. 
\end{proof}
\end{lem}

In Theorem \ref{AIP} we show that $B_t^n$ can be approximated by a sum of $Z_j$ which, by Lemma \ref{Zcent}, are i.i.d.\ centred random variables. With the aim of proving a central limit theorem, we now show that they also have finite second moments. 
\begin{lem}\label{ZfinVar}
Suppose that $\beta>1$ and $\Eb[\eta_0^2]<\infty$ then $\Eb[Z_j^2]<\infty$ for $j\geq 2$.
\begin{proof}
Since $\{Z_j\}_{j\geq 2}$ are i.i.d.\ under $\Pb$ we have that $\mathrm{Var}_\Pb(Z_j)=\mathrm{Var}_\Pb(Z_2) $ for all $j\geq 2$. By properties of regenerations times $Y_{\kappa_2}\geq Y_{\kappa_1} $ and $S_{\kappa_2}\geq S_{\kappa_1} $ almost surely therefore we have that 
\begin{flalign}\label{sigsq}
\mathrm{Var}_\Pb(Z_2)  \leq \Eb\left[\left(Y_{\kappa_2}-Y_{\kappa_1}\right)^2\right] + \nu^2\Eb\left[\left(S_{\kappa_2}-S_{\kappa_1}\right)^2\right]. 
\end{flalign}
For the second term we have
\begin{flalign}
 \Eb\left[\left(S_{\kappa_2}-S_{\kappa_1}\right)^2\right] & = \Eb\left[\!\left(\sum_{x=Y_{\kappa_1}}^{Y_{\kappa_2}-1}\sum_{i=1}^{L(x,\infty)}\eta_{x,i}\right)^2\right] \notag \\
 & = \Eb\left[\sum_{x=Y_{\kappa_1}}^{Y_{\kappa_2}-1}\!\left(\sum_{i=1}^{L(x,\infty)}\eta_{x,i}\right)^2\right]+\Eb\left[\sum_{x=Y_{\kappa_1}}^{Y_{\kappa_2}-1}\sum_{y=Y_{\kappa_1}}^{Y_{\kappa_2}-1}\ind_{\{x\neq y\}}\!\left(\sum_{i=1}^{L(x,\infty)}\eta_{x,i}\right)\!\!\!\left(\sum_{j=1}^{L(y,\infty)}\eta_{y,j}\right)\!\right].\label{Sbnd}
\end{flalign}
By conditioning on $Y$ we have that the holding times at separate vertices are independent therefore the second term in this expression can be written as
\begin{flalign*}
 & \Eb\left[\sum_{x=Y_{\kappa_1}}^{Y_{\kappa_2}-1}\sum_{y=Y_{\kappa_1}}^{Y_{\kappa_2}-1}\ind_{\{x\neq y\}}\Eb\left[\sum_{i=1}^{L(x,\infty)}\eta_{x,i}\Big| Y\right]\Eb\left[\sum_{j=1}^{L(y,\infty)}\eta_{y,j}\Big|Y\right]\right] \\
 & \qqqqqqquad =  \Eb\left[\sum_{x=Y_{\kappa_1}}^{Y_{\kappa_2}-1}\sum_{y=Y_{\kappa_1}}^{Y_{\kappa_2}-1}\ind_{\{x\neq y\}}\left(\sum_{i=1}^{L(x,\infty)}\Eb[\eta_{x,i}]\right)\!\!\!\left(\sum_{j=1}^{L(y,\infty)}\Eb[\eta_{y,j}]\right)\right] \\
 &  \qqqqqqquad = \Eb\left[\sum_{x=Y_{\kappa_1}}^{Y_{\kappa_2}-1}\sum_{y=Y_{\kappa_1}}^{Y_{\kappa_2}-1}\ind_{\{x\neq y\}}\Eb[\eta_0]^2L(x,\infty)L(y,\infty)\right] \\
 &  \qqqqqqquad \leq \Eb[\eta_0]^2 \Eb\left[\left(\sum_{x=Y_{\kappa_1}}^{Y_{\kappa_2}-1}L(x,\infty)\right)^2\right] \\
 & \qqqqqqquad  = \Eb[\eta_0]^2 \Eb[(\kappa_2-\kappa_1)^2].
\end{flalign*}
By (\ref{expMom}), the time between regenerations $\kappa_2-\kappa_1$ has exponential moments therefore $\Eb[(\kappa_2-\kappa_1)^2]<\infty$. Furthermore, since $Y$ moves in discrete time and has jumps of length $1$ 
\begin{flalign}\label{Yrho}
\Eb\left[\left(Y_{\kappa_2}-Y_{\kappa_1}\right)^2\right] \leq \Eb\left[\left(\kappa_2-\kappa_1\right)^2\right]<\infty. 
\end{flalign}
Combining (\ref{Yrho}) with (\ref{sigsq}) and (\ref{Sbnd}), in order to show that $\mathrm{Var}_{\Pb}(Z_2)<\infty$ it remains to show that 
\[ \Eb\left[\sum_{x=Y_{\kappa_1}}^{Y_{\kappa_2}-1}\!\left(\sum_{i=1}^{L(x,\infty)}\eta_{x,i}\right)^2\right]<\infty.\]
Conditioning on $Y$ this expectation is equal to
\begin{flalign*}
 \Eb\left[\sum_{x=Y_{\kappa_1}}^{Y_{\kappa_2}-1}\sum_{i,j=1}^{L(x,\infty)}\Eb\left[\eta_{x,i}\eta_{x,j}\Big|Y\right]\right] \; \leq \;  \Eb\left[\sum_{x=Y_{\kappa_1}}^{Y_{\kappa_2}-1}L(x,\infty)^2\Eb\left[\eta_{x,1}^2\right]\right] \; \leq \; \Eb[\eta_0^2]\Eb[(\kappa_2-\kappa_1)^2]
\end{flalign*}
which is finite by the assumptions of the theorem and equation (\ref{Yrho}).
\end{proof}
\end{lem}

We now conclude the proof of the annealed functional central limit theorem by showing that $B_t^n$ can be suitably approximated by a sum of $Z_j$. 
\begin{thm}\label{AIP}
 Suppose that $\beta>1$ and $\Eb[\eta_0^2]<\infty$ then there exists $\varsigma^2\in(0,\infty)$ such that
 \[B_t^n=\frac{X_{nt} -nt\nu}{\varsigma \sqrt{n}}\] converges in $\Pb$-distribution on $D(\Rb^+,\Rb)$ endowed with the Skorohod metric to a standard Brownian motion.  
 \begin{proof}
By Lemmas \ref{Zcent} and \ref{ZfinVar}
 \[\Sigma_m:= \sum_{j=2}^{m}Z_j=\left(X_{S_{\kappa_m}}-S_{\kappa_m}\nu\right) -\left(X_{S_{\kappa_1}}-S_{\kappa_1}\nu\right)\] 
 for $m\geq 2$ is a sum of i.i.d.\ centred random variables with finite second moment. 
  
Write $m_t:=\sup\{j\geq 0:S_{\kappa_j}\leq t\}$ to be the number of regenerations by time $t>0$ then
\begin{flalign*}
 &\sup_{t \in [0,T]}\left|B_t^n-\frac{\Sigma_{m_{tn}}}{\varsigma \sqrt{n}}\right| \; \leq \; \frac{X_{S_{\kappa_1}}+S_{\kappa_1}+|\min_kY_k|}{\varsigma\sqrt{n}}+\sup_{j=1,...,m_{Tn}}\frac{X_{S_{\kappa_{j+1}}}-X_{S_{\kappa_j}}+\left(S_{\kappa_{j+1}}-S_{\kappa_j}\right)\nu}{\varsigma\sqrt{n}}.
\end{flalign*}
The random variables $X_{S_{\kappa_1}},S_{\kappa_1}$ and $|\min_kY_k|$ are all almost surely finite therefore the first fraction converges to $0$ $\Pb$-a.s. For $\varepsilon>0$, by a union bound and Markov's inequality
\begin{flalign*}
 &\Pb\left(\sup_{j=1,...,m_{Tn}}\frac{X_{S_{\kappa_{j+1}}}-X_{S_{\kappa_j}}}{\sqrt{n}}>\varepsilon\right)\\
 &\qquad \leq \Pb\left(m_{Tn}>2Tn\Eb[\eta_0]^{-1}\right)+C_Tn\Pb\left(X_{S_{\kappa_2}}-X_{S_{\kappa_1}}>\varepsilon\sqrt{n}\right) \\
 &\qquad \leq \Pb\left(m_{Tn}>2Tn\Eb[\eta_0]^{-1}\right)+C_{T,\varepsilon}\Eb\left[(X_{S_{\kappa_2}}-X_{S_{\kappa_1}})^2\ind_{\{X_{S_{\kappa_2}}-X_{S_{\kappa_1}}\geq \varepsilon\sqrt{n}\}}\right].
\end{flalign*}
By Proposition \ref{Sconv}, since $S^{-1}_t\geq m_t$, we have that $\Pb\left(m_{Tn}>2Tn\Eb[\eta_0]^{-1}\right)\rightarrow 0$ as $\nin$. By (\ref{Yrho}) $\Eb\left[(X_{S_{\kappa_2}}-X_{S_{\kappa_1}})^2\right]<\infty$ therefore by dominated convergence 
\[\Eb\left[(X_{S_{\kappa_2}}-X_{S_{\kappa_1}})^2\ind_{\{X_{S_{\kappa_2}}-X_{S_{\kappa_1}}\geq \varepsilon\sqrt{n}\}}\right]\rightarrow 0\] 
as $\nin$. Similarly, by Lemma \ref{ZfinVar}, $\Eb\left[(S_{\kappa_2}-S_{\kappa_1})^2\right]<\infty$ hence we have that \[\Pb\left(\sup_{j=1,...,m_{tn}}\frac{S_{\kappa_{j+1}}-S_{\kappa_j}}{\sqrt{n}}>\varepsilon\right) \rightarrow 0\]
as $\nin$, and the supremum distance between $(B_t^n)_{t \in [0,T]}$ and $(\Sigma_{m_{tn}}/\varsigma\sqrt{n})_{t \in [0,T]}$ converges to $0$ in $\Pb$-probability. It therefore suffices to prove an invariance principle for $\Sigma_{m_{tn}}$.

For $s\in \Rb^+$ let $\Sigma_s$ denote the linear interpolation of $\Sigma_m$ then by Donsker's invariance principle we have that $(\Sigma_{t n}/\sqrt{n})_{t \in [0,T]}$ converges in distribution to a scaled Brownian motion. 

By the law of large numbers we have that $\kappa_n/n$ converges $\Pb$-a.s.\ to $\Eb[\kappa_2-\kappa_1]$ as $\nin$. Therefore, by \cite[Corollary 13.6.4]{wh}, we have that $m_{tn}/n$ converges $\Pb$-a.s.\ on $D([0,\infty),\Rb)$ with the Skorohod metric to the deterministic process $R_t:=(\Eb[\eta_0]\Eb[\kappa_2-\kappa_1])^{-1}t$. 

By \cite[Theorem 13.2.1]{wh} it follows that the sequence $(\Sigma_{m_{tn}}/\sqrt{n})_{t \in [0,T]}$ converges to the same limit as $(\Sigma_{R_tn}/\sqrt{n})_{t \in [0,T]}$ which is a scaled Brownian motion. In particular, choosing 
\begin{flalign}\label{varsig}\varsigma^2=\frac{\Eb[Z_2^2]}{\Eb[\eta_0]\Eb[\kappa_2-\kappa_1]}\end{flalign}
we have that $B_t^n$ converges to a standard Brownian motion.
 \end{proof}

\end{thm}

\section{Quenched central limit theorem}\label{QuCLT}
In this section we prove a quenched central limit theorem for the randomly trapped random walk under a $2+\varepsilon$ moment condition. We do this by adapting the method used in \cite{go} and first proving a quenched CLT for the first hitting time of $n$. Write $\tau_n:=\inf\{t\geq 0:X_t=n\}$ and, for $\omega$ fixed, $\Hc^\omega(n):=\Et^\omega[\tau_n]$. Let $\zeta_k:=\tau_{k+1}-\tau_k$ be the time taken between hitting $k$ and $k+1$ for the first time then the elements of $(\zeta_k)_{k\geq 1}$ are independent under $\Pt^\omega$ and \[\tau_n=\sum_{k=0}^{n-1}\zeta_k.\]
\begin{lem}\label{HitCLT}
 Suppose that $\beta>1$ and $\Eb[\eta_0^2]<\infty$, then for $\Pr$-a.e.\ $\omega$ we have that \[\Pt^\omega\left(\frac{\tau_n-\Hc^\omega(n)}{\sigma\sqrt{n}} < x\right) \rightarrow \Phi(x)\] 
 uniformly in $x$ as $\nin$, where $\sigma^2=\Er[\mathrm{Var}_{\Pt^\omega}(\tau_1)]$ and \[\Phi(x)=\int_{-\infty}^x\frac{e^{-\frac{u^2}{2}}}{\sqrt{2\pi}} \d u.\]
 \begin{proof}
  
  By definition of $\tau_n, \Hc^\omega(n)$ and $\zeta_k$ \[\frac{\tau_n-\Hc^\omega(n)}{\sigma\sqrt{n}}=\frac{\sum_{k=0}^{n-1}(\zeta_k-\Et^{\omega}[\zeta_k])}{\sigma\sqrt{n}}.\]
  It therefore suffices to show that Lindeberg's conditions (see \cite[Theorem 3.4.5]{du}) hold:
  \begin{enumerate}
   \item for $\Pr$-a.e.\ $\omega$, as $\nin$ \[\sum_{k=0}^{n-1}\Et^\omega\left[\left(\frac{\zeta_k-\Et^\omega[\zeta_k]}{\sigma\sqrt{n}}\right)^2\right]\rightarrow 1;\]
   \item for $\Pr$-a.e.\ $\omega$, $\forall \varepsilon>0$ as $\nin$ \[\sum_{k=0}^{n-1}\Et^\omega\left[\left(\frac{\zeta_k-\Et^\omega[\zeta_k]}{\sigma\sqrt{n}}\right)^2\ind_{\{|\zeta_k-\Et^\omega[\zeta_k]|>\varepsilon\sqrt{n}\}}\right]\rightarrow 0.\]
  \end{enumerate}

Recall, from the remark prior to Proposition \ref{Sconv}, that $\theta$ is the shift map which is ergodic by \cite[Lemma 2.1]{cewa}. For the first condition we have that $\zeta_k-\Et^\omega[\zeta_k]=\theta^k(\zeta_0-\Et^\omega[\zeta_0])$. These random variables are identically distributed under $\Pr$ with $\Er\left[\Et^\omega[(\zeta_0-\Et^\omega[\zeta_0])^2]\right]<\infty$ therefore
    \begin{flalign*}
   \sum_{k=0}^{n-1}\frac{\mathrm{Var}_\omega\left(\zeta_k\right)}{\sigma^2n} & = \sum_{k=0}^{n-1}\frac{\mathrm{Var}_{\theta^k\omega}\left(\zeta_0\right)}{\sigma^2n}
  \end{flalign*}
which converges to $\Er\left[\mathrm{Var}_\omega(\zeta_0)\right]\sigma^{-2}=1$ for $\Pr$-a.e.\ $\omega$ by Birkhoff's ergodic theorem. 

For the second condition write $U_K^\omega(\cdot):=\Et^\omega[(\cdot-\Et^\omega[\cdot])^2\ind_{\{|\cdot-\Et^\omega[\cdot]|>K\}}]$ then for all $\varepsilon>0$ there exists $\exists N_{\varepsilon,K}\in \Nb$ such that $\varepsilon\sqrt{n}>K$ for all $n\geq N_{\varepsilon,K}$. Therefore, for $n$ large 
\begin{flalign*}
 \sum_{k=0}^{n-1}\Et^\omega\left[\left(\frac{\zeta_k-\Et^\omega[\zeta_k]}{\sigma\sqrt{n}}\right)^2\ind_{\{|\zeta_k-\Et^\omega[\zeta_k]|>\varepsilon\sqrt{n}\}}\right] 
 \; \leq \; \sum_{k=0}^{n-1}\frac{U_K^\omega(\zeta_k)}{\sigma^2n} 
  \;=\; \sum_{k=0}^{n-1}\frac{U_K^{\theta^k\omega}(\zeta_0)}{\sigma^2n}. 
\end{flalign*}
By Birkhoff's ergodic theorem, for $\Pr$-a.e.\ $\omega$, this converges to 
\[\frac{\Er\left[U_K^{\omega}(\zeta_0)\right]}{\sigma^2}=\frac{\Er\left[\Et^\omega\left[\left(\zeta_0-\Et^\omega[\zeta_0]\right)^2\ind_{\{|\zeta_0-\Et^\omega[\zeta_0]|>K\}}\right]\right]}{\sigma^2}\] which converges to $0$ as $K\rightarrow \infty$ by dominated convergence.
%
 \end{proof}

\end{lem}

Write $\tau^Y_n:=\inf\{m\geq 0: \; Y_m=n\}$ to be the first hitting time of level $n$ by the underlying walk. The following lemma describes the probability that the underlying walk moves back $k$ levels before moving forward $n$. The result is the classical Gambler's ruin therefore we omit the proof.

\begin{lem}\label{returnProb}
For integers $k<0<n$
\begin{flalign*}
\Pt_0(\tau_k^Y<\tau_n^Y) = \frac{\beta^n-1}{\beta^{n-k}-1}.
\end{flalign*}
\end{lem}

By Lemma \ref{HitCLT} we have a central limit theorem for the first hitting time of vertex $n$. The environment dependent centring $\Hc^\omega(n)$ can be written as the sum of $n$ identically distributed random variables $\Et^\omega[\zeta_k]$. These are not independent however; they are only locally dependent. Recall that $\eta_{k,i}$ is the $i^{\text{th}}$ holding time at vertex $k$ hence $\Er[\zeta_0]=\Er[\eta_{0,0}](\beta+1)/(\beta-1)$ then write 
\[\tilde{\Hc}^\omega(n):=\sum_{k=0}^{n-1}\frac{\beta+1}{\beta-1}\Et^\omega[\eta_{k,0}].\]
We now show that $\Hc^\omega$ and $\tilde{\Hc}^\omega$ do not differ too much and therefore Lemma \ref{HitCLT} also holds with $\Hc^\omega$ replaced by $\tilde{\Hc}^\omega$. Under $\Pr$, the function $\tilde{\Hc}^\omega(n)$ is a sum of i.i.d.\ random variables with non-zero variance unless $\Et^\omega[\eta_0]$ is constant. We thus have a central limit theorem for $\tilde{\Hc}^\omega$ (and therefore $\Hc^\omega$), which shows that the environment dependent centring is necessary. Notice that this is the first point at which we introduce the extra $2+\varepsilon$ moment condition however we do require the condition later in Lemma \ref{condJ} as well.
\begin{lem}\label{HHtil}
 Suppose $\beta>1$ and that $\Er\left[\Et^\omega[\eta_0]^{2+\varepsilon}\right]<\infty$ for some $\varepsilon>0$, then for $\Pr$-a.e.\ $\omega$ \[\left|\frac{\tilde{\Hc}^\omega(n)-\Hc^\omega(n)}{\sqrt{n}}\right|\rightarrow 0.\]
 \begin{proof}
  Recall that $L(k,m)$ denotes the local time of $Y$ at vertex $k$ by time $m$ and that the trapping times $\eta_{k,j}$ do not depend on the underlying walk, then 
  \begin{flalign*}
  \Hc^\omega(n) \; = \; \Et^\omega\left[\sum_{k=-\infty}^{n-1}\sum_{j=1}^{L(k,\tau_n^Y)}\eta_{k,j}\right] \; = \; \sum_{k=-\infty}^{n-1}\Et_0[L(k,\tau_n^Y)]\Et^\omega[\eta_{k,0}].
  \end{flalign*}
We need to determine the expected local times at sites up to reaching level $n$. That is, $\Et_0[L(k,\tau_n^Y)]=\Pt_0(\tau_k^Y<\tau_n^Y)\Et_k[L(k,\tau_n^Y)]$ (by the strong Markov property).

Let $(\tau^Y_n)^+:=\inf\{m> 0: \; Y_m=n\}$ be the first return time to level $n$ by the underlying walk. By Lemma \ref{returnProb}, the number of visits to $k$ before reaching $n$ for a walk started at $k$ is geometrically distributed with escape probability 
\[\Pt_k(\tau_n^Y<(\tau_k^Y)^+)=\frac{\beta}{1+\beta}\left(1-\Pt_0(\tau_{-1}^Y<\tau_{n-k-1}^Y)\right)=\frac{\beta^{n-k}(\beta-1)}{(\beta^{n-k}-1)(\beta+1)}.\]
Therefore, 
\[\Et_k[L(k,\tau_n^Y)]=\frac{(\beta^{n-k}-1)(\beta+1)}{\beta^{n-k}(\beta-1)}\]
and 
\[\Et_0[L(k,\tau_n^Y)]= \begin{cases}  \frac{(\beta^n-1)(\beta+1)}{\beta^{n-k}(\beta-1)}    & \text{if } k<0 \\ \frac{(\beta^{n-k}-1)(\beta+1)}{\beta^{n-k}(\beta-1)} & \text{if } k\geq 0.                 \end{cases}\]

For fixed $k<0$, $\Et_0[L(k,\tau_n^Y)]$ is increasing in $n$ and converges to $\beta^k(\beta+1)/(\beta-1)$. In particular,
\begin{flalign*}
 0  \leq \sum_{k=-\infty}^{-1}\Et_0[L(k,\tau_n^Y)]\Et^\omega[\eta_{k,0}] 
  \leq C \sum_{k=1}^\infty \beta^{-k}\Et^\omega[\eta_{-k,0}]
\end{flalign*}
which is finite for $\Pr$-a.e.\ $\omega$ therefore $n^{-1/2}\sum_{k=-\infty}^{-1}\Et_0[L(k,\tau_n^Y)]\Et^\omega[\eta_{k,0}]$ converges to $0$ for $\Pr$-a.e.\ $\omega$.

For $k\geq 0$ fixed, $\Et_0[L(k,\tau_n^Y)]$ is increasing in $n$ and converges to $(\beta+1)/(\beta-1)$. In particular,
\begin{flalign*}
 0 & \leq \sum_{k=0}^\infty \left(\frac{\beta+1}{\beta-1}-\Et_0[L(k,\tau_n^Y)]\right)\Et^\omega[\eta_{k,0}] \\
  & = \frac{\beta+1}{\beta-1}\sum_{k=0}^{n-1} \beta^{-(n-k)}\Et^\omega[\eta_{k,0}] \\
  & = \frac{\beta+1}{\beta-1}\left(\sum_{k=0}^{n-\left\lfloor2\frac{\log(n)}{\log(\beta)}\right\rfloor} \beta^{-(n-k)}\Et^\omega[\eta_{k,0}]+ \sum_{k=n-\left\lfloor2\frac{\log(n)}{\log(\beta)}\right\rfloor+1}^{n-1} \beta^{-(n-k)}\Et^\omega[\eta_{k,0}]\right) \\
  & \leq C\left(\sum_{k=0}^{n-1}\beta^{-2\frac{\log(n)}{\log(\beta)}}\Et^\omega[\eta_{k,0}]+\sum_{k=n-\left\lfloor2\frac{\log(n)}{\log(\beta)}\right\rfloor+1}^{n-1}\Et^\omega[\eta_{k,0}]\right) \\
  & = C\left(\sum_{k=0}^{n-1}\frac{\Et^\omega[\eta_{k,0}]}{n^2}+\sum_{k=n-\left\lfloor2\frac{\log(n)}{\log(\beta)}\right\rfloor+1}^{n-1}\Et^\omega[\eta_{k,0}]\right). 
\end{flalign*}
The first term converges to $0$ for $\Pr$-a.e.\ $\omega$ by the strong law of large numbers. For the second term we have that, for $\delta,\epsilon>0$, by Markov's inequality 
 \begin{flalign*}
  \Pr\left(\sum_{k=n-\left\lfloor2\frac{\log(n)}{\log(\beta)}\right\rfloor+1}^{n-1}\Et^\omega[\eta_{k,0}]>\epsilon\sqrt{n}\right) & \leq 2\frac{\log(n)}{\log(\beta)}\Pr\left(\Et^\omega[\eta_0]>\frac{\epsilon\sqrt{n}\log(\beta)}{2\log(n)}\right) \\
  & = 2\frac{\log(n)}{\log(\beta)}\Pr\left(\Et^\omega[\eta_0]^{2+\delta}>\frac{C_\epsilon n^{1+\delta/2}}{\log(n)^{2+\delta}}\right) \\
 & \leq \frac{C\log(n)^{3+\delta}}{n^{1+\delta/2}}
 \end{flalign*}
since we can choose $\delta>0$ sufficiently small such that $\Er\left[\Et^\omega[\eta_0]^{2+\delta}\right]<\infty$. By Borel-Cantelli we then have that \[\sum_{k=n-\left\lfloor2\frac{\log(n)}{\log(\beta)}\right\rfloor+1}^{n-1}\frac{\Et^\omega[\eta_{k,0}]}{\sqrt{n}}\] converges to $0$ for $\Pr$-a.e.\ $\omega$.
\end{proof}
\end{lem}

We now prove a technical lemma that allows us to control the difference between $\tilde{\Hc}^\omega$ and its expected value under $\Pr$ which is important in proving the quenched CLT for the walk.
\begin{lem}\label{condJ}
 Let \[\Jc(n):=\sum_{j=0}^{n-1}\left(\Et^\omega[\eta_{k,0}]-\Eb[\eta_{k,0}]\right)\]
 and \[\Jc^*(n):=\max_{m\leq n}\Jc(m).\]
\begin{enumerate}
 \item\label{opc} Suppose $\Er\left[\Et^\omega[\eta_0]^{2}\right]<\infty$, then for any $c>0$, $\Jc(n)n^{-\frac{1+c}{2}} \rightarrow 0$ for $\Pr$-a.e.\ $\omega$;
 \item\label{tpd} Suppose $\Er\left[\Et^\omega[\eta_0]^{2+\varepsilon}\right]<\infty$ for some $\varepsilon>0$, then for $\delta>0$ sufficiently small and some constant $C$ \[\Er\left[|\Jc^*(n)|^{2+2\delta}\right]^{\frac{1}{2+2\delta}}\leq Cn^{1/2}.\]
\end{enumerate}
\begin{proof}
 By Theorem 17 of Chapter IX.3 in \cite{pe}, if $Z_n$ are i.i.d.\ centred random variables, $a_n$ is an increasing, diverging sequence and 
 \begin{enumerate}
  \item\label{ChebBC} $\sum_{n=1}^\infty \Pr(|Z_1|\geq a_n)<\infty$;
  \item\label{inttest} $\sum_{k=n}^\infty a_k^{-2}=O\left(\frac{n}{a_n^2}\right)$;
  \item\label{knbnd} $a_k/a_n \leq Ck/n$ for all $k\geq n$
 \end{enumerate}
then $\sum_{k=1}^n \frac{Z_k}{a_n}$ converges to $0$, $\Pr$-a.s. 

Write $Z_n:=\Et^\omega[\eta_{n,0}]-\Eb[\eta_{n,0}]$ then $Z_n$ are i.i.d.\ and centred under $\Pr$; moreover, the sequence $a_n=n^{\frac{1+c}{2}}$ is increasing and diverges. By Chebyshev we have that
\[\sum_{n=1}^\infty \Pr(|Z_1|\geq a_n) \leq \sum_{n=1}^\infty \mathrm{Var}_{\Pr}(\Et^\omega[\eta_0])n^{-(1+c)} <\infty \]
which gives condition \ref{ChebBC}. Since $n/a_n^2=n^{-c}$, an integral test gives condition \ref{inttest}. For $k\geq n$ we have that $a_k/a_n =(k/n)^{\frac{1+c}{2}} \leq k/n$ so long as $c\leq 1$ which gives \ref{knbnd}. We therefore have that for any $c>0$, $\Jc(n)n^{-\frac{1+c}{2}} \rightarrow 0$ for $\Pr$-a.e.\ $\omega$ hence the first statement holds.

The process $\Jc(m)$ is a martingale therefore by the $L^p$-maximal inequality we have that 
\[\Er\left[\max_{m\leq n}|\Jc(m)|^{2+2\delta}\right]\leq \left(\frac{2+2\delta}{1+2\delta}\right)^{2+2\delta}\Er\left[|\Jc(n)|^{2+2\delta}\right].\]
It therefore suffices to show that 
\[\Er\left[\left(\sum_{j=0}^{n-1}\frac{\Et^\omega[\eta_{k,0}]-\Eb[\eta_{k,0}]}{\sqrt{n}}\right)^{2+2\delta}\right]\]
is bounded above. By the Marcinkiewicz-Zygmund inequality \cite[Theorem 5]{mazy} we have that 
\[\Er\left[\left(\sum_{j=0}^{n-1}\frac{\Et^\omega[\eta_{k,0}]-\Eb[\eta_{k,0}]}{\sqrt{n}}\right)^{2+2\delta}\right]\leq C\Er\left[\left(\sum_{j=0}^{n-1}\left(\frac{\Et^\omega[\eta_{k,0}]-\Eb[\eta_{k,0}]}{\sqrt{n}}\right)^2\right)^{1+\delta}\right]\]
which is bounded above by
\[C\Er\left[\sum_{j=0}^{n-1}\frac{\left(\Et^\omega[\eta_{k,0}]-\Eb[\eta_{k,0}]\right)^{2+2\delta}}{n}\right]=C\Er\left[\left(\Et^\omega[\eta_0]-\Eb[\eta_0]\right)^{2+2\delta}\right]\]
using Jensen's inequality. Using that $\Er\left[\Et^\omega[\eta_0]^{2+\varepsilon}\right]<\infty$ for some $\varepsilon>0$, it then follows that for $\delta>0$ sufficiently small and some constant $C$ 
\[\Er\left[|\Jc^*(n)|^{2+2\delta}\right]^{\frac{1}{2+2\delta}}\leq Cn^{1/2}.\]
\end{proof}
\end{lem}

We now prove the main result of the section which is a quenched central limit theorem for the randomly trapped random walk. Recall from Corollary \ref{LinX} that $\nu=(\beta-1)\left((\beta+1)\Eb[\eta_0]\right)^{-1}$ is the $\Pb$-a.s.\ limit of $X_n/n$ and from Lemma \ref{HitCLT} that $\varsigma^2=\Er[\mathrm{Var}_\omega(\tau_1)]$ is the variance in the quenched CLT for the first hitting times $\tau_n$; write $\vartheta:=\varsigma \nu^{3/2}$ and
\[\Gc^\omega(t):=\nu t-\nu\sum_{k=0}^{\lfloor \nu t-1\rfloor}\frac{\beta+1}{\beta-1}(\Et^\omega[\eta_{k,0}]-\Eb[\eta_{k,0}]). \]

Notice that, since $\nu t$ is the centring of $X_t$ in the annealed CLT (Theorem \ref{AIP}) and $\Gc^\omega(t)-\nu t$ is a sum of i.i.d.\ random variables under the environment law, whenever these random variables have non-zero variance we obtain a central limit theorem for $\Gc^\omega(t)-\nu t$ with respect to $\Pr$ and therefore proving Theorem \ref{QCLT} also shows that there is no quenched CLT for $X_t$ with a deterministic centring. 
\begin{thm}\label{QCLT}
 Suppose $\beta>1$, $\Eb[\eta_0^2]<\infty$ and for some $\varepsilon>0$ we have that $\Er[\Et^\omega[\eta_0]^{2+\varepsilon}]<\infty$, then for $\Pr$-a.e.\ $\omega$ we have that \[\Pt^\omega\left(\frac{X_t-\Gc^\omega(t)}{\vartheta\sqrt{t}} \leq x\right) \rightarrow \Phi(x)\]
 uniformly in $x$ as $\nin$.
 \begin{proof}
  Let $\varrho_t :=\sup\{|X_s|: \; s\leq t\}$ be the furthest point reached by $X$ up to time $t$; then $\tau_{\varrho_t}\leq t<\tau_{\varrho_t+1}$ and $X_t\leq X_{\tau_{\varrho_t}}=\varrho_t$. We then have that $|X_t-\varrho_t|=|X_t-X_{\tau_{\varrho_t}}| \leq \sup_{s\geq \tau_{\varrho_t}}X_{\tau_{\varrho_t}}-X_s$. Since $X_{nt}/n$ converges $\Pb$-a.s.\ to $\nu t$ uniformly over $t$ we can choose a constant $C_T$ such that for $n$ sufficiently large $X_t\leq C_Tn$ for all $t\leq nT$. Write \[A_n:=\bigcap_{k=1}^{\lfloor C_Tn+1\rfloor}\left\{\inf\{Y_m:m\geq \tau_k^Y\} \geq k-C\log(n)\right\}\] 
  to be the event that the walk never backtracks distance $C\log(n)$ up to reaching vertex $\lceil C_Tn\rceil$. By Lemma \ref{returnProb} we then have that
  \[\Pt(A_n^c) \; \leq \; C_Tn\Pt(\tau_{-\lfloor C\log(n)\rfloor}<\infty) \; \leq \; C_Tn\beta^{-C\log(n)} \; = \; C_Tn^{1-C\log(\beta)}.\]
  Therefore, choosing $C$ such that $C\log(\beta)> 2$, by Borel-Cantelli we have that there exists only finitely many $n$ such that the walk backtracks distance $C\log(n)$ up to time $nT$. In particular, on this event $|X_t-\varrho_t|n^{-1/2}\leq C\log(n)n^{-1/2}$ which converges deterministically to $0$ uniformly over $t\leq T$. It therefore suffices to show that for $\Pr$-a.e.\ $\omega$ \[\lim_{t\rightarrow \infty}\Pt^\omega\left(\frac{\varrho_t-\Gc^\omega(t)}{\vartheta\sqrt{t}}\leq x\right)=\Phi(x).\]
  
  By monotonicity we have that $\{\varrho_t\leq m\}=\{\tau_{m+1}>t\}$. Writing $\Ic^\omega(t):=\lfloor x\vartheta\sqrt{t}+\Gc^\omega(t)+1\rfloor$ it then follows that
  \begin{flalign*}
   \Pt^\omega\left(\frac{\varrho_t-\Gc^\omega(t)}{\vartheta\sqrt{t}}<x\right) \; = \; \Pt^\omega\left(\tau_{\Ic^\omega(t) }>t\right) \; = \; \Pt^\omega\left(
   \frac{\tau_{\Ic^\omega(t)}-\Hc^\omega(\Ic^\omega(t))}{\sigma\sqrt{\Ic^\omega(t)}}>\frac{t-\Hc^\omega(\Ic^\omega(t))}{\sigma\sqrt{t}}\cdot \sqrt{\frac{t}{\Ic^\omega(t)}}
   \right).
  \end{flalign*}
The sequence $\Ic^\omega(t)$ is increasing in $t$ and diverges; in particular, by the law of large numbers $t/\Ic^\omega(t)$ converges to $\nu^{-1}$ for $\Pr$-a.e.\ $\omega$. The result then follows from Lemma \ref{HitCLT} if $\Hc^\omega(\Ic^\omega(t))=t+\sigma \nu^{1/2} x\sqrt{t}+o_t$, where $o_t/\sqrt{t}$ converges to $0$ for $\Pr$-a.e.\ $\omega$.

Since $\Ic^\omega(t)$ diverges, by Lemma \ref{HHtil} it suffices to show that $\tilde{\Hc}^\omega(\Ic^\omega(t))=t+\sigma \nu^{1/2} x\sqrt{t}+o_t$. By definition of $\tilde{\Hc}^\omega$ and $\Ic^\omega(n)$ we have that there exists some $O_1:=O_1(\omega,t,x)$ such that $|O_1| \leq  \nu^{-1}$ and
\begin{flalign*}
 \tilde{\Hc}^\omega(\Ic^\omega(t)) & = \nu^{-1}\Ic^\omega(t)  + \sum_{k=0}^{\Ic^\omega(t)-1} \frac{\beta+1}{\beta-1}(\Et^\omega[\eta_{k,0}]-\Eb[\eta_{k,0}]) \\
 & = t+\sigma \nu^{1/2} x \sqrt{t} - \sum_{k=0}^{\lfloor \nu t\rfloor-1} \frac{\beta+1}{\beta-1}(\Et^\omega[\eta_{k,0}]-\Eb[\eta_{k,0}])  +\sum_{k=0}^{\Ic^\omega(t)-1} \frac{\beta+1}{\beta-1}(\Et^\omega[\eta_{k,0}]-\Eb[\eta_{k,0}]) +O_1. 
\end{flalign*}
Moreover, for some $O_2:=O_2(\omega,t,x)$ satisfying $|O_2| \leq 3$, we have that
\[\Ic^\omega(t)-\lfloor\nu t\rfloor =\vartheta x\sqrt{t}+\Eb[\eta_0]^{-1}\sum_{k=0}^{\lfloor\nu t-1\rfloor}\left(\Et^\omega[\eta_{k,0}]-\Eb[\eta_{k,0}]\right) + O_2.\]

By part \ref{opc} of Lemma \ref{condJ} we have that $(\Ic^\omega(t)-\lfloor\nu t\rfloor)t^{-\frac{1+c}{2}}$ converges to $0$ for $\Pr$-a.e.\ $\omega$ and any $c>0$. In order to show that $\tilde{\Hc}^\omega(\Ic^\omega(t))=t+\sigma \nu^{1/2} x\sqrt{t}+o_t$ it now suffices to show that for all $c>0$ suitably small 
\[\Rc^\omega(n,c):=n^{-1/2}\max_{m\leq n^{\frac{1+c}{2}}}\left|\sum_{k=n}^{n+m}\left(\Et^\omega[\eta_{k,0}]-\Eb[\eta_{k,0}]\right)\right|\]
converges to $0$ for $\Pr$-a.e.\ $\omega$. 

Suppose that $\Rc^\omega(n^2,2c)$ converges to $0$ for all $c>0$ suitably small and $\Pr$-a.e.\ $\omega$ then for $i =1,...,2n$ (that is, $n^2< n^2+i <(n+1)^2$) we have that
\begin{flalign*}
 \Rc^\omega(n^2+i,c) & = (n^2+i)^{-1/2}\!\!\!\max_{m\leq (n^2+i)^{\frac{1+c}{2}}}\left|\sum_{k=n^2+i}^{n^2+i+m}\!\left(\Et^\omega[\eta_{k,0}]-\Eb[\eta_{k,0}]\right)\right| \\
 & \leq (n^2+i)^{-1/2}\!\!\!\max_{m\leq (n^2+i)^{\frac{1+c}{2}}}\left(\left|\sum_{k=n^2}^{n^2+i+m}\!\left(\Et^\omega[\eta_{k,0}]-\Eb[\eta_{k,0}]\right)\right| +\left|\sum_{k=n^2}^{n^2+i-1}\!\left(\Et^\omega[\eta_{k,0}]-\Eb[\eta_{k,0}]\right)\right|\right).
\end{flalign*}
Since $i+m<n^{1+2c}$ for $n$ suitably large we then have that $\Rc^\omega(n^2+i,c)\leq 2\Rc^\omega(n^2,2c)$ for all $i =1,...,2n$ thus it suffices to show that $\Rc^\omega(n^2,2c)$ converges to $0$ for all $c>0$ suitably small and $\Pr$-a.e.\ $\omega$. For $\epsilon>0$, by Markov's inequality
\begin{flalign*}
 \Pr\left(\Rc^\omega(n^2,2c)>\epsilon\right) & \leq \Er\left[\Rc^\omega(n^2,2c)^{2+2\delta}\right]\epsilon^{-(2+2\delta)} \\
 & = C_\epsilon\Er\left[\left(n^{-1}\max_{m\leq n^{1+2c}}\left|\sum_{k=n^2}^{n^2+m}\left(\Et^\omega[\eta_{k,0}]-\Eb[\eta_{k,0}]\right)\right|\right)^{2+2\delta}\right]  \\
 & \leq C_\epsilon n^{(1+2c)(1+\delta)-2(1+\delta)}
\end{flalign*}
by part \ref{tpd} of Lemma \ref{condJ}. Choosing $c<\delta/(2+2\delta)$ gives us that \[\sum_{n=1}^\infty\Pr\left(\Rc^\omega(n^2,2c)>\epsilon\right)<\infty \]
therefore by Borel-Cantelli we have the desired result.

 \end{proof}
\end{thm}

\section{Subcritical Galton-Watson trees}\label{Sub}
A subcritical Galton-Watson tree conditioned to survive consists of a semi-infinite path with GW-trees as leaves. The walk on the tree does not deviate too far from this path and therefore behaves like a randomly trapped random walk on $\Zb$ with holding times distributed as excursion times in GW-trees. Biased walks on subcritical GW-trees are, therefore, a natural example of the randomly trapped random walk. In this section we prove the required conditions for a biased random walk on a subcritical GW-tree conditioned to survive to satisfy the convergence results of the previous sections. That is, we prove a speed result, an Einstein relation, an annealed functional CLT and a quenched CLT with an environment dependent centring. To begin, we set up the model and show that it suffices to consider a randomly trapped random walk on $\Zb$ with trapping times distributed as excursion times in trees. 

Let $f(s)=\sum_{k=0}^\infty p_ks^k$ denote the generating function of a GW-process with mean $\mu\in(0,1)$ and variance $\sigma^2 <\infty$. To avoid the trivial case in which no traps form, we also assume that there exists $k\geq 2$ such that $p_k>0$. We denote $\xi$ as a random variable with the offspring law $\Pr(\xi=k)=p_k$ and $\xi^*$ a random variable with the size-biased law given by the probabilities $\Pr(\xi^*=k)=kp_k\mu^{-1}$. Let $Z_n$ denote the $n\th$ generation size of the $f$-GW-process with this law. That is, started from a single progenitor, each individual independently has $k$ offspring with probability $p_k$. Such a process gives rise to a random rooted tree $\Tc$ where individuals in the process are represented by vertices (with the unique progenitor as the root $\rho$) and undirected edges connect individuals with their offspring. More generally, we will use $Z_n^T$ to denote the $n\th$ generation size of a tree $T$. It has been shown in \cite{ke} that there is a well defined probability measure $\Pr$ over $f$-GW trees conditioned to survive which arises as the limit as $\nin$ of probability measures over GW-trees conditioned to survive up to generation $n$. It can be seen (e.g.\ \cite{ja}) that the tree can be constructed by attaching i.i.d.\ finite trees to a single infinite path $\Yc:=(\rho_0=\rho,\rho_1,...)$. More specifically, starting with a single special vertex $\rho_0$, at each generation let every normal vertex give birth to normal vertices according to independent copies of the original offspring distribution and every special vertex give birth to vertices according to independent copies of the size-biased distribution, one of which is chosen uniformly at random to be special (and denoted $\rho_k$ in the $k\th$ generation). We will use $\Tc$ to denote an $f$-GW-tree and $\Tc^*$ an $f$-GW-tree conditioned to survive. We refer to $\Yc$ as the backbone of the tree and the finite connected components $\Tc^*\setminus\{(\rho_i,\rho_{i+1}); \; i\geq 0\}$ appended to $\Yc$ as branches. We denote by $\Tc^{*-}_x$, the branch rooted at $x \in \Yc$; that is, the descendants of $x$ which are not in the descendant tree of the unique child of $x$ on the backbone. 

\begin{figure}[H]
\centering
 \includegraphics[scale=0.8]{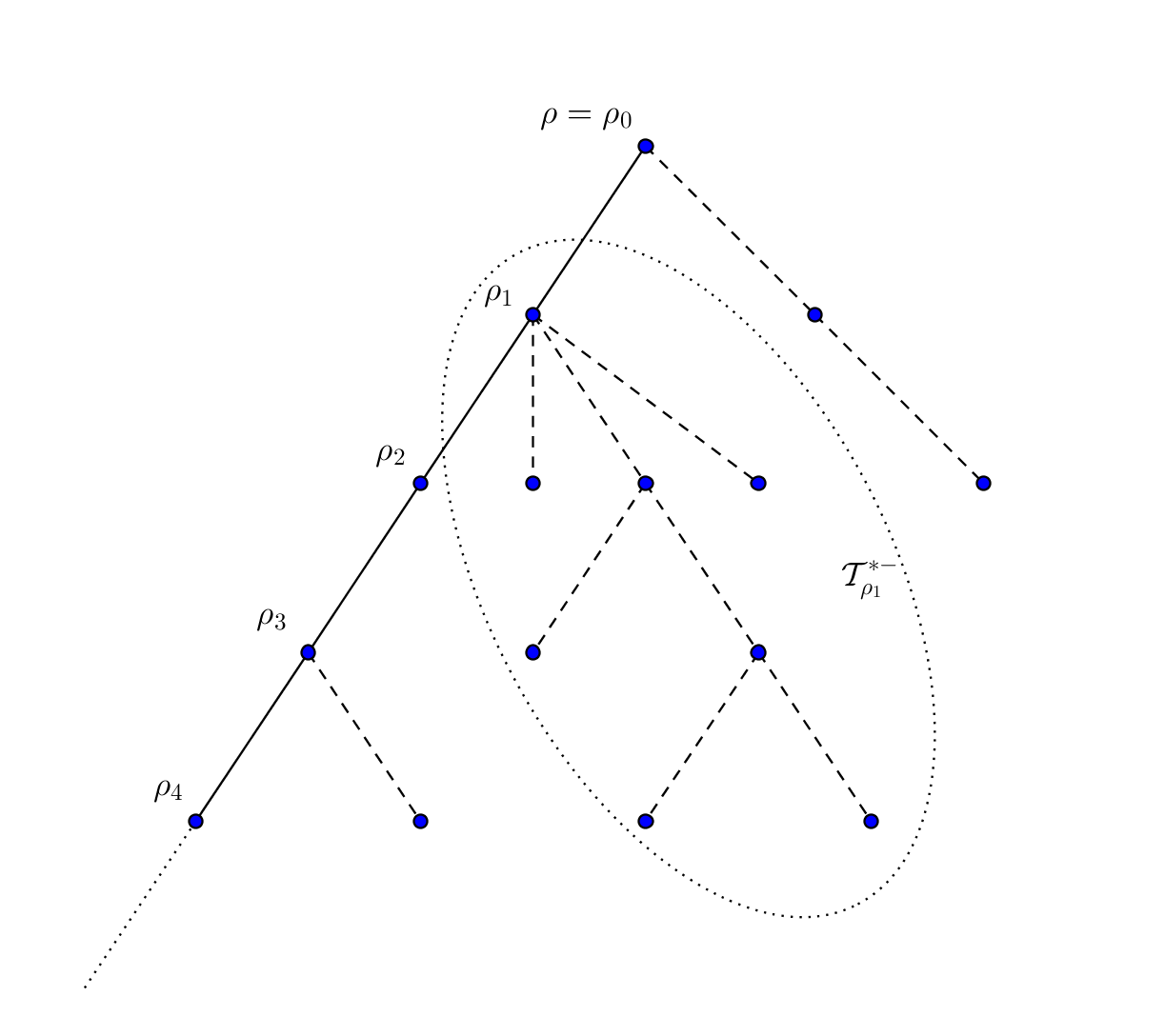} 
\caption{A sample subcritical GW-tree conditioned to survive $\Tc^{*-}$ with the backbone $\Yc$ represented by solid lines and the buds and traps connected by dashed lines.}\label{treediag}
\end{figure}

A $\beta$-biased random walk on a fixed, rooted tree $T$ is a random walk $(X_n)_{n \geq 0}$ on $T$ which is $\beta$-times more likely to make a transition to a given child of the current vertex than the parent (which are the only options). More specifically, let $\rho$ denote the root of the $T$, $\overleftarrow{x}$ the parent of $x \in T$, $c(x)$ the set of children of $x$ and $d_x:=|c(x)|$ the number of children, then the random walk $X$ is the Markov chain started from $X_0=z$ defined by the transition probabilities 
\[\Pt^T_z(X_{n+1}=y|X_n=x)=\begin{cases} \frac{1}{1+\beta d_x}, & \text{if } y=\overleftarrow{x}, \\  \frac{\beta}{1+\beta d_x}, & \text{if } y \in c(x), \; x \neq \rho, \\ \frac{1}{
d_\rho}, & \text{if } y \in c(x), \; x =\rho, \\ 0, & \text{otherwise.} \\ \end{cases} \]

We use $\Pb_\rho(\cdot)=\int \Pt^{\Tc^*}_\rho(\cdot)\Pr(\text{d}\Tc^*)$ for the annealed law obtained by averaging the quenched law $\Pt^{\Tc^*}_\rho$ over a law $\Pr$ on $f$-GW-trees conditioned to survive (with a fixed root $\rho$). For the remainder of the section we shall assume that $\beta>1$ so that the walk is $\Pr$-a.s.\ transient.

In \cite{bowd} it is shown that either a strong bias or heavy tails of the offspring law can slow the walk into a sub-ballistic phase. More specifically, when $\beta\geq \mu^{-1}$ the walk spends a large amount of time at the deep parts of the largest traps and when $\sigma^2=\infty$ the walk takes a large number of excursions into the branches. These slowing effects cause the walk to become sub-ballistic; that is, $|X_n|n^{-1}$ converges $\Pb$-a.s.\ to $0$ in either case. Further to this, an exponent $\gamma$ has been determined such that the walk escapes polynomially with this exponent.

Let $|X_n|$ denote the graph distance between the walk at time $n$ and the root of the tree. We show that if $\beta \in (1,\mu^{-1})$ and $\Er[\xi^2]<\infty$ then $|X_n|n^{-1}$ converges $\Pb$-a.s.\ to
\[\nu=\frac{\mu(\beta-1)(1-\beta\mu)}{\mu(\beta+1)(1-\beta\mu)+2\beta(\sigma^2-\mu(1-\mu))}.\]
We then prove that if $\Er[\xi^3]<\infty$ and $\beta \in (1,\mu^{-1/2})$ then
\[B_t^n:=\frac{|X_{nt}|-nt\nu}{\varsigma \sqrt{n}}\]
converges in distribution under $\Pb$ to a Brownian motion. Following this we show that if $\beta \in (1,\mu^{-1/2})$ and $\Er[\xi^{3+\varepsilon}]<\infty$ for some $\varepsilon>0$ then for $\Pr$-a.e.\ tree $\Tc^*$,
\[\frac{|X_{n}|-\Gc^{\Tc^*}(n)}{\vartheta \sqrt{n}}\]
converges in distribution to a standard Gaussian under $\Pt^{\Tc^*}$ where 
\[\Gc^{\Tc^*}(t)=\nu t - \nu\sum_{k=1}^{\lfloor \nu t\rfloor}\frac{\beta+1}{\beta-1}\left(\Et^{\Tc^*}[\tilde{\eta}_{\rho_k,0}]-\Eb[\eta_0]\right)\]
is an environment dependent centring. Moreover, under these conditions, $\Gc^{\Tc^*}(n)/\sqrt{n}$ converges in distribution under $\Pr$ to a Gaussian random variable with strictly positive variance. This confirms the need for an environment dependent centring. We also show that the conditions $\beta^2\mu< 1$ and $\Er[\xi^3]<\infty$ are necessary for the previous results.

In Section \ref{SRTRW} we couple the walk on the tree with a randomly trapped random walk in such a way that our main results of this section can be deduced from corresponding results for the RTRW. We then derive moment bounds on the excursion times in random trees that allow us to apply the results of the previous sections. Asymptotics for the generation size of a GW-tree will be important throughout the section; before discussing the walk in detail we give several useful properties. An important relation we shall use throughout is that the sequence $\Pr(Z_n>0)\mu^{-n}$ is decreasing and converges, which follows from \cite[Theorem B]{lypepeot} since $\mu \in (0,1)$ and $\sigma^2<\infty$, but can also be traced back further to \cite{heseve} and \cite{ya}. In particular, we shall write by $c_\mu$, the limiting constant:
\begin{flalign}\label{cmu}
\limn \frac{\Pr(Z_n>0)}{\mu^{n}}=c_\mu.
\end{flalign}

Lemma \ref{gensize} shows bounds on the expected moments of the generation sizes. A simple extension shows that these upper bounds are tight.
\begin{lem}\label{gensize}
 Let $Z_n$ denote the $n^{th}$ generation size of an $f$-GW-process with offspring distribution $\xi$ and mean $\mu \in (0,1)$. 
 \begin{enumerate}
  \item\label{EZn} $\Er[Z_n]=\mu^n$.
  \item\label{ZZ} If $\Er[\xi^2]<\infty$ and $m\geq n$ then $\Er[Z_nZ_m]\leq C\mu^m$ for some constant $C$.
  \item\label{ZZZ} If $\Er[\xi^3]<\infty$ and $l\geq m\geq n$ then $\Er[Z_nZ_mZ_l]\leq C\mu^l$ for some constant $C$.
 \end{enumerate}
\begin{proof}
If $\Pr(\xi< 2)=1$ then $Z_n$ only takes values $0$ and $1$. Moreover, the number of generations until extinction is geometrically distributed with termination probability $1-\mu$ therefore the result follows. We therefore assume that $\Pr(\xi\geq 2)>0$ which implies that $f''(1)>0$.

Let $f_n$ denote the generating function of $Z_n$ then statement \ref{EZn} follows from
\begin{flalign}\label{EZ1}
\Er\left[Z_n\right] \; = \; \sum_{j=1}^\infty j\Pr\left(Z_n=j\right)\; = \;f'_n(1)\; = \;f'_{n-1}(1)f'(n)\; = \;f'(1)^n\; = \;\mu^n.
\end{flalign}

For the second moment we have that $\Er[Z_n^2]=f''_n(1)+\Er[Z_n]$ where
 \begin{flalign*}
  f''_{n+1}(1) \; = \; \left(f\left(f_n(s)\right)'\right)'\Big|_{s=1}  \; = \; f''(1)\mu^{2n}+\mu f''_n(1).
  \end{flalign*}
  Applying this recursively we see that
  \begin{flalign*}
  f''_{n+1}(1) \; = \; f''(1)\sum_{k=0}^n\mu^{n+k}  \; = \; f''(1)\mu^n\frac{1-\mu^{n+1}}{1-\mu}
 \end{flalign*}
therefore 
\begin{flalign}\label{Z2}
\Er[Z_n^2]= c_1\mu^n+c_2\mu^{2n}
\end{flalign}
whenever $f''(1)<\infty$ which follows from $\Er[\xi^2]<\infty$. For $m>n$, by stationarity of GW-processes and (\ref{EZ1}) we have that
\[\Er[Z_m|Z_n=k] \;= \; \Er[Z_{m-n}|Z_0=k] \; = \; k\Er[Z_{m-n}] \; = \; k\mu^{m-n}.\]
Therefore, statement \ref{ZZ} follows by 
\begin{flalign}\label{ZmZnk}
 \Er[Z_nZ_m] \; = \; \sum_{k=1}^\infty k \Pr(Z_n=k)\Er[Z_m|Z_n=k]  \; = \; \sum_{k=1}^\infty k^2 \Pr(Z_n=k)\mu^{m-n}   \; = \; \mu^{m-n}\Er[Z_n^2] \; \leq \; C\mu^m.
\end{flalign}

When $\Er[\xi^3]<\infty$ we have that $f''(1),f'''(1)<\infty$ therefore differentiating $f''_n(s)$ and evaluating at $s=1$ gives us that
\begin{flalign*}
 f'''_{n+1}(1)& = 3f''(1)f'_n(1)f''_n(1)+f'_n(1)^3f'''(1)+f'(1)f'''_n(1)  \\
 & = \left(f'''(1)-\frac{3f''(1)^2}{\mu(1-\mu)}\right)\mu^{3n}+\frac{3f''(1)^2}{\mu(1-\mu)} \mu^{2n}+\mu f'''_n(1).
\end{flalign*}
Iterating gives us that
\[f'''_{n+1}(1)=\left(f'''(1)-\frac{3f''(1)^2}{\mu(1-\mu)}\right)\frac{1-\mu^{2n}}{1-\mu}\mu^{3n}+\frac{3f''(1)^2}{\mu(1-\mu)} \cdot\frac{1-\mu^n}{1-\mu}\mu^{n}+\mu^nf'''(1)\] 
which proves that, for some $c_1>0$, 
\begin{flalign}\label{Z3}
\Er[Z_n^3]=c_1\mu^n+c_2\mu^{2n}+c_3\mu^{3n}+c_4\mu^{5n}
\end{flalign} 

If $l\geq m\geq n$ and $\Er[\xi^3]<\infty$ then for any $j,k\geq 1$ (where $j=k$ if $m=n$)
\[\Er[Z_nZ_mZ_l|Z_m=j,Z_n=k] \; = \; kj\Er[Z_l|Z_m=j] \; = \; kj^2\Er[Z_{l-m}] \; = \; kj^2\mu^{l-m}\]
by stationarity and (\ref{EZ1}) therefore
\begin{flalign}
 \Er[Z_nZ_mZ_l]  & = \sum_{k=1}^\infty  \Pr(Z_n=k)\sum_{j=1}^\infty  \Pr(Z_m=j|Z_n=k) \Er[Z_nZ_mZ_l|Z_m=j,Z_n=k] \notag \\
 & =\mu^{l-m} \sum_{k=1}^\infty k \Pr(Z_n=k)\sum_{j=1}^\infty j^2 \Pr(Z_m=j|Z_n=k) \notag \\
 & = \mu^{l-m}\sum_{k=1}^\infty k \Pr(Z_n=k)\Er[Z_m^2|Z_n=k]. \label{Znml}
 \end{flalign}
For $j\geq 1$ let $Z_n^{(j)}$ be independent copies of $Z_n$ then by convexity and (\ref{Z2})
\begin{flalign}\label{Zm2Zn}
\Er[Z_m^2|Z_n=k] \;= \; k^2\Er\left[\left(\sum_{j=1}^k\frac{Z_{m-n}^{(j)}}{k}\right)^2\right] \;\leq \; k^2\Er[Z_{m-n}^2] \; \leq Ck^2\mu^{m-n}
\end{flalign}
therefore, combining this with (\ref{Znml}) we have that
\[ \Er[Z_nZ_mZ_l]  \leq C\mu^{l-n} \Er[Z_{n}^3]  \leq C\mu^l,\]
where the final inequality follows by (\ref{Z3}), which gives statement \ref{ZZZ}.

\end{proof}
\end{lem}

\subsection{Describing the walk on the subcritical tree as a randomly trapped random walk}\label{SRTRW}
We now wish to construct an almost equivalent model which allows us to consider the walk on the GW-tree in our randomly trapped random walk framework. To begin, we show that the walk never deviates too far from the backbone. For a fixed tree $T$ with root $\rho$ we write $\Hc(T):=\max_{x \in T}d(\rho,x)$ to be the height of the tree. Let $\tilde{X}_n$ be the projection of $X_n$ onto $\Yc$; that is, $\tilde{X}_n$ is the unique vertex on $\Yc$ which satisfies $d(X_n,\tilde{X}_n)=\min_{y \in \Yc}d(X_n,y)$. 
\begin{lem}\label{proBack}
Suppose $\mu\in(0,1), \beta\geq 1$ and $\sigma^2<\infty$ then $\Pb$-a.s.\
 \[\sup_{n\geq 2}\sup_{m\leq nT}\frac{d(X_m,\tilde{X}_m)}{\log(n)}<\infty.\]
 \begin{proof}
  The distance $d(X_m,\tilde{X}_m)$ between the walk and the backbone is at most the height of the largest branch seen up to time $nT$ therefore, since the walk can have seen at most $M$ branches by time $M$, by a union bound we have that for $C>0$
  \begin{flalign*}
   \Pb\left(\sup_{m\leq nT}d(X_m,\tilde{X}_m)> C\log(n)\right) &  \leq \lceil nT\rceil \Pr\left(\Hc(\Tc_\rho^{*-})> C\log(n)\right) \\
   & \leq \lceil nT\rceil \Er[\xi^*]\Pr\left(\Hc(\Tc)> C\log(n)-1\right). 
  \end{flalign*}
By properties of the size-biased distribution we have that $\Er[\xi^*]=\mu^{-1}\Er[\xi^2]<\infty$, since $\sigma^2<\infty$. Therefore, by (\ref{cmu}),
\begin{flalign*}
 \lceil nT\rceil \Er[\xi^*]\Pr\left(\Hc(\Tc)> C\log(n)\right) & \leq C_T n \mu^{C\log(n)}.
\end{flalign*}
We can therefore choose $C$ sufficiently large so that \[\Pb\left(\sup_{m\leq nT}d(X_m,\tilde{X}_m)> C\log(n)\right) \leq C_Tn^{-2}\]
thus the result follows by Borel-Cantelli.
 \end{proof}
\end{lem}

For $x \in T$ recall that $|x|:=d(\rho,x)$ then $|\tilde{X}_n|$ has the same distribution as a randomly trapped random walk on $\Nb$ with holding times distributed as excursion times in trees. The traps formed at different vertices are independent and identically distributed except at $\rho$ since the root does not have an ancestor and therefore the transition probabilities from $\rho$ differ from those at other vertices on the backbone. We now show that we can extend from $\Nb$ to $\Zb$ with i.i.d.\ traps. 

We begin by constructing the holding times of the randomly trapped random walk via a sequence of i.i.d.\ trees. Start with an initial vertex $\rho$ and a unique ancestor $\oR$. Attach $\xi^*-1$ offspring to $\rho$ where $\xi^*$ is size-biased as above. Note that this could result in zero offspring of $\rho$ in which case the tree ends with only vertices $\rho,\oR$. Otherwise, attach independent $f$-GW trees to the offspring of $\rho$. This creates a tree $\oT$ which has the distribution of a branch with an additional vertex connected to the root. 

\begin{figure}[H]
\centering
 \includegraphics[scale=0.8]{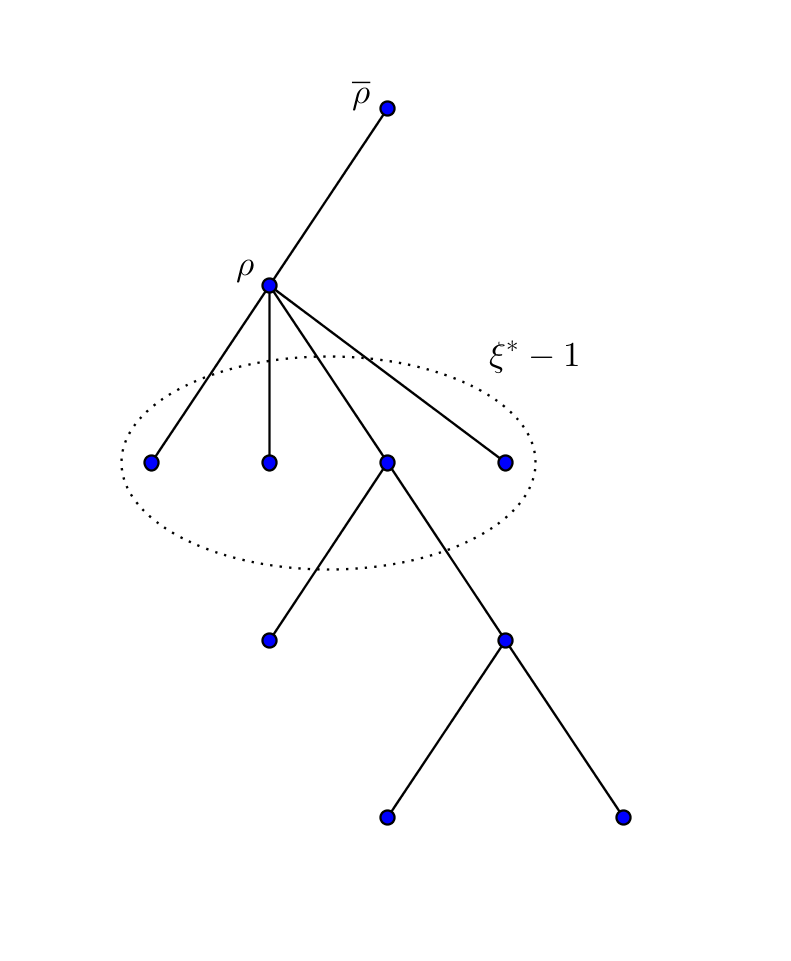} 
\caption{A tree $\oT$ with fixed vertices $\rho,\oR$ and $\xi^*-1$ independent $f$-GW-trees attached to $\rho$.}
\end{figure}

Recall that $\overleftarrow{x}$ denotes the parent of $x \in T$ and $c(x)$ the set of children of $x$. Consider a walk $(W_n)_{n\geq 0}$ on $\oT$ with transition probabilities 
\[\Pt^{\oT}(W_{n+1}=y|W_n=x)=\begin{cases} 1, & \text{if } x=y=\oR, \\ \frac{\beta+1}{\beta (|c(x)|+1)+1}, & \text{if } x=\rho, y=\oR, \\  \frac{\beta}{\beta (|c(x)|+1)+1}, & \text{if } x=\rho, y \in c(x), \\ \frac{\beta}{|c(x)|\beta+1}, & \text{if } x \notin\{\rho,\oR\}, y \in c(x),  \\ \frac{1}{|c(x)|\beta+1}, & \text{if } x \notin\{\rho,\oR\}, y =\overleftarrow{x},  \\ 0, & \text{otherwise.} \\ \end{cases} \] 
An excursion in $\oT$ started from $\rho$ until absorption in $\oR$ has the same distribution as the time taken to move between backbone vertices of $\Tc^*$ (except at the root of $\Tc^*$). Let $\omega=(\oT_x)_{x \in \Zb}$ denote a sequence of independent trees with this law. For $\omega$ fixed let $(\eta_{x,i})_{x \in \Zb,i\geq 0}$ be independent with 
\[\Pt^\omega\left(\eta_{x,i}=k\right)=\Pt^{\oT_x}_\rho\left(\min\{n>0:W_n=\oR\}=k\right)\]
where $\rho,\oR$ are the vertices in $\oT_x$ corresponding with the construction.

Recall that for a discrete time process $W$ we let 
\[L^W(x,n):=\sum_{k=0}^n\ind_{\{W_k=x\}}\]
denote its local time. Let $r(0)=0$ and for $k=1,2,...$ define $r(k):=\inf\{n>r(k-1):\tilde{X}_n\neq \tilde{X}_{r(k-1)}$ to be the time of the $k^{th}$ movement of $\tilde{X}$. The walk $\tilde{Y}_n:=\tilde{X}_{r(n)}$ is then a $\beta$-biased walk on $\Yc$ reflected at $\rho$. Moreover, for $x \in \Zb$ and $i=1,2,...$ we can write 
\[\tilde{\eta}_{x,i}:=r(k+1)-r(k) \quad \text{where} \quad k=\min\{j:L^{\tilde{Y}}(x,j)=i\}\]
to be the holding time of $\tilde{X}$ at vertex $x$ on the $i\th$ visit.

Let  $\hat{Y}_n$ be a simple, $\beta$-biased random walk on $\Zb$, then define 
\[A_n:=\sum_{k=1}^n\ind_{\{\hat{Y}_k,\hat{Y}_{k-1}\geq 0\}} \qquad \text{and} \qquad A_n^{-1}:=\sup\{m\geq 0:A_m\leq n\}.\]
The process $\hat{Y}_{A_n^{-1}}$ is equal in distribution to $|\tilde{Y}_n|$ therefore, without loss of generality, we may couple $\tilde{Y}$ to $\hat{Y}$ in the construction of $X$ so that $|\tilde{Y}_n|=\hat{Y}_{A_n^{-1}}$ without changing the distribution of $X$. 

Let \[\hat{S}_n:=\sum_{x \in \Zb}\sum_{i=1}^{L^{\hat{Y}}(x,n-1)}\hat{\eta}_{x,i} \qquad \text{where} \qquad \hat{\eta}_{x,i}=\begin{cases} \eta_{x,i}, & \text{if } x \leq 0, \\ \tilde{\eta}_{\rho_x,i}, & \text{if } x>0. \end{cases} \]
Write $\hat{S}_n^{-1}:=\inf\{k\geq 0:\hat{S}_k>n\}$ then $\hat{X}_n:=\hat{Y}_{\hat{S}^{-1}_n}$ is a randomly trapped random walk coupled to $\tilde{X}$ with trapping times equal in distribution to $(\eta_{x,i})_{x \in \Zb,i\geq 0}$. The following lemma shows that $\tilde{X}$ and $\hat{X}$ never deviate too far.

\begin{lem}\label{Back}
  If $\beta\mu<1$ and $\beta>1$ then for any $\delta>0$ we have that \[\frac{\sup_{m\leq nT}||\tilde{X}_m|-\hat{X}_m|}{n^\delta}\]
 converges $\Pb$-a.s.\ to $0$.
 \begin{proof}
Using that $|\tilde{X}|,\hat{X}$ are discrete time processes with jump size $1$, by the coupling of the two process we have that 
\begin{flalign}\label{TilHatBnd}
\sup_{m\leq nT}||\tilde{X}_m|-\hat{X}_m|
\leq \sum_{k=0}^\infty\ind_{\{\tilde{Y}_k=\rho_0\}}\tilde{\eta}_{\rho_0,L^{\tilde{Y}}(\rho_0,k)} +\sum_{x\leq 0}\sum_{k=0}^\infty\ind_{\{\hat{Y}_k=x\}}\eta_{x,L^{\hat{Y}}(x,k)} 
\end{flalign}
which is independent of $n$. That is, the supremum distance between the two processes is at most the total time spent by the two processes where the holding times differ.
 
Recall that $\tau_x^+:=\inf\{m>0: X_m=x\}$ is the first return time to $x$ for a $\beta$ biased random walk. For a fixed tree $T$ rooted at $\rho$ with $n\th$ generation size $Z_n$ where $Z_1>0$ it is classical (e.g. \cite[Chapter 2]{lype}) that 
\begin{flalign}\label{return}
\Et^T_\rho[\tau^+_\rho]=2\sum_{n\geq 1}\frac{Z_n\beta^{n-1}}{Z_1}.
\end{flalign}
The holding times $\eta_{x,k}$ are distributed as first hitting times of vertex $\oR$ by the walk $W_n$ on $\oT$ started from vertex $\rho$. The transition probabilities only differ from those of a $\beta$-biased walk at $\rho$ and $\oR$ but $W_n$ is more likely to move to $\oR$ than the $\beta$-biased walk. It follows that the holding times are stochastically dominated by first return times to $\oR$ of a $\beta$-biased random walk on a copy of $\oT$ started from $\oR$. This means that 
\[\Eb[\eta_{x,1}] \leq 2\sum_{n\geq 1}\Er\left[Z_{n}^{\oT}\right]\beta^{n-1} =2\left(1+\beta\Er[\xi^*-1]\sum_{n=0}^\infty \Er\left[Z_n^{\Tc}\right]\beta^n\right) \]
where $\Tc$ is an $f$-GW-tree and we've used that for $n\geq 2$
\begin{flalign*}
 \Er\left[Z_n^{\oT}\right]  = \Er\left[\Er\left[Z_n^{\oT}|Z_2^\oT=\xi^*-1\right]\right] = \Er\left[(\xi^*-1)\Er\left[Z_{n-2}^{\Tc}\right]\right]=\Er\left[\xi^*-1\right]\Er\left[Z_{n-2}^{\Tc}\right]
\end{flalign*}
which follows from independence and stationarity of GW-trees. By statement \ref{EZn} of Lemma \ref{gensize} we have that $\Er\left[Z_n^\Tc\right]=\mu^n$ therefore
\[\Eb[\eta_{x,1}] \leq 2+\frac{2\beta}{1-\beta\mu}\Er[\xi^*-1] \]
which is finite since $\beta\mu<1$ and $\sigma^2<\infty$. Similarly, the holding times $\tilde{\eta}_{0,k}$ are stochastically dominated by excursion times of a $\beta$-biased random walk on $\oT$ started from $\oR$ therefore we have that $\Eb[\tilde{\eta}_{\rho_0,1}]<\infty$.

For a $\beta$-biased walk on $\Zb$ the last hitting time of $0$ has exponential moments since it must occur before the first regeneration time and by (\ref{expMom}) the time of the first regeneration has exponential moments. Therefore, using the above bound we can conclude that 
\[\Eb\left[\sum_{k=0}^\infty\ind_{\{\tilde{Y}_k=\rho_0\}}\tilde{\eta}_{\rho_0,L^{\tilde{Y}}(\rho_0,k)} +\sum_{x\leq 0}\sum_{k=0}^\infty\ind_{\{\hat{Y}_k=x\}}\eta_{x,L^{\hat{Y}}(x,k)}\right]<\infty \]
which is sufficient to prove the result by (\ref{TilHatBnd}).
 \end{proof}
\end{lem}

As a result of Lemmas \ref{proBack} and \ref{Back} we can consider the randomly trapped random walk model in order to prove results for a biased random walk on a subcritical GW-tree conditioned to survive. In the remainder of the section we derive conditions for the tree which yield the moment bounds required in Theorems  \ref{AIP} and \ref{QCLT} and Corollary \ref{LinX}.

\subsection{The speed of the walk}
We now wish to prove a bound on the expected holding time for the randomly trapped random walk. Let $\eta_k:=\eta_{\hat{Y}_k,L(\hat{Y}_k,k)}$ be the $k\th$ holding time of the randomly trapped random walk. Recall that $\eta_0$ is distributed as the first hitting time of $\oR$ by the walk $W_n$ on the tree $\oT$ described at the beginning of the section, and that $\oT$ is a random tree rooted at $\rho$ with a single ancestor $\oR$ of the root, $\xi^*-1$ buds attached as children of $\rho$ (where $\xi^*$ has the size biased law) and independent $f$-GW-trees attached to the buds.

The quantity $\eta_0$ is distributed as the first hitting time of $\oR$ in the random tree $\oT$ by the walk $W$ started from $\rho$. Let 
 \begin{flalign}\label{Nlab}
 N=\sum_{k=1}^{\tau_{\oR}}\ind_{\{W_k=\rho\}}
 \end{flalign}
be the number of return times to the root $\rho$ before reaching its unique ancestor $\oR$. That is, $N$ is the number of excursions to the trees attached to $\rho$ before the walk reaches $\oR$. Let $\tau_x^{(0)}:=0$ then for $k=1,...,N$ write $\tau_x^{(k)}:=\min\{n> \tau_x^{(k-1)}:W_n=x\}$ to be the hitting times of $x$ and $\zeta_{k}:=\tau_{\rho}^{(k)}-\tau_{\rho}^{(k-1)}$ the duration of the $k^{\text{th}}$ excursion. We then have that 
\begin{flalign}\label{eta1}
 \eta_0=1+\sum_{k=1}^N\zeta_{k}.
\end{flalign}

\begin{lem}\label{spd}
 Suppose $\beta\mu<1$, $\sigma^2<\infty$ and $\beta\geq1$, then \[\Eb[\eta_0]=\frac{\mu(\beta+1)(1-\beta\mu)+2\beta(\sigma^2-\mu(1-\mu))}{\mu(\beta+1)(1-\beta\mu)}.\]
 Moreover, if $\sigma^2=\infty$ or $\beta\geq \mu^{-1}$ then $\Eb[\eta_0]=\infty$.
 \begin{proof}
By (\ref{eta1}) we have that
\begin{flalign*}
 \Eb[\eta_0] = 1+\Eb\left[\sum_{k=1}^N \Eb[\zeta_k|N]\right] = 1+\Eb[N]\Er\left[\Et^{\Tc}_\rho[\tau_\rho^+]|Z_1=1\right]
\end{flalign*}
where $Z_1$ denotes the first generation size of an $f$-GW-tree $\Tc$. The number of excursions $N$ is geometrically distributed under $\Pt^{\oT_0}$ with termination probability $1-p_{ex}$ where
\[p_{ex}:=\Pt^{\oT_0}(W_1\neq \oR)=\frac{\beta(\xi^*-1)}{\beta\xi^*+1}.\]
It therefore follows that
\begin{flalign*}
 \Eb[N] & = \Eb\left[\frac{\beta(\xi^*-1)}{\beta+1}\right]  = \frac{\beta}{\beta+1}\left(\sum_{k=1}^\infty \frac{k^2p_k}{\mu}-1\right)  = \frac{\beta(\sigma^2-\mu(1-\mu))}{(\beta+1)\mu}.
\end{flalign*}

Using the formula (\ref{return}) for the expected time spent in a fixed tree and statement \ref{EZn} of Lemma \ref{gensize} for the expected size of the $k\th$ generation we have that
\begin{flalign*}
 \Er\left[\Et^{\Tc}_\rho[\tau_\rho^+]|Z_1=1\right] = \Er\left[2\sum_{k\geq 1}\frac{Z_k\beta^{k-1}}{Z_1}|Z_1=1\right]=2\sum_{k\geq 1}\Er[Z_k|Z_1=1]\beta^{k-1} = 2\sum_{k\geq 1}(\beta\mu)^{k-1}.
\end{flalign*}
If $\beta<\mu^{-1}$ then this is equal to $2/(1-\beta\mu)$; otherwise, the sum does not converge. It follows that 
\[\Eb[\eta_0]=\frac{\mu(\beta+1)(1-\beta\mu)+2\beta(\sigma^2-\mu(1-\mu))}{\mu(\beta+1)(1-\beta\mu)}.\]
 \end{proof}
\end{lem}

By Lemma \ref{proBack} and \ref{Back} we have that $\sup_{m \leq nT}||X_m|-\hat{X}_m|n^{-1}$ converges to $0$ $\Pb$-a.s.\ therefore the following result follows from Corollary \ref{LinX} and Lemma \ref{spd} since $\hat{X}$ has the distribution of a randomly trapped random walk with trapping times equal in distribution to $(\eta_{x,i})_{x \in \Zb,i\geq 0}$.

\begin{cly}\label{SubTreeSpd}
 Suppose $\beta\mu<1,\; \sigma^2<\infty$ and $\beta>1$, then $|X_{n}|/n$ converges $\Pb$-a.s.\ to $\nu_\beta$ where \[ \nu_\beta := \frac{\mu(\beta-1)(1-\beta\mu)}{\mu(\beta+1)(1-\beta\mu)+2\beta(\sigma^2-\mu(1-\mu))}.\]
\end{cly}

We now extend the Einstein relation for the randomly trapped random walk to the walk on the GW-tree. This is a non-trivial extension because, in the tree model, the bias affects the trapping times and the unbiased walk is significantly influenced by the restriction to the half line. For this reason we observe convergence to a reflected Brownian motion and cannot simply apply Corollary \ref{EIN}.
\begin{lem}\label{SubEIN}
 Suppose $\mu<1$ and $\sigma^2<\infty$. The unbiased $(\beta=1)$ walk $X_{\lfloor nt \rfloor}n^{-1/2}$ converges in $\Pb$-distribution on $D([0,\infty),\Rb)$ endowed with the Skorohod metric to $|B_t|$ where $B_t$ is a scaled Brownian motion with variance $\Upsilon=\Eb[\eta_0]^{-1}$. Moreover, 
\[\lim_{\beta\rightarrow 1^+}\frac{\nu_\beta}{\beta-1}=\frac{\Upsilon}{2}\]
where $\nu_\beta$ is the speed calculated in Corollary \ref{SubTreeSpd} for the $\beta$-biased walk.
\begin{proof}
Recall that $\hat{X}_n$ is an randomly trapped random walk on $\Zb$ which, by assumption and Lemma \ref{spd}, is unbiased and has finite expected holding times. By \cite[Theorem 2.9]{arcacero}, for $\Pr$-a.e. $\omega$, the rescaled process $\hat{X}_{nt}n^{-1/2}$ converges in $\Pt^\omega$ distribution to a scaled Brownian motion $B$ with variance $\Eb[\eta_0]^{-1}$. 

The scaled local time at the origin $n^{-1}L^{\hat{Y}}(0,n-1)$ converges $\Pt$-a.s.\ to $0$. Moreover, the holding times $(\eta_{0,i})_{i\geq1}$ are i.i.d.\ under $\Pt^\omega$ therefore, by the law of large numbers, $\sum_{i=1}^nn^{-1}\eta_{0,i}$ converges $\Pt^\omega$-a.s.\ to $\Et^\omega[\eta_{0,1}]$ for $\Pr$-a.e.\ $\omega$. The same holds for $(\tilde{\eta}_{0,i})_{i\geq1}$ therefore the scaled sums 
\[\sum_{i=1}^{L^{\hat{Y}}(0,n-1)}\frac{\eta_{0,i}}{n} \qquad \text{and} \qquad \sum_{i=1}^{L^{\hat{Y}}(0,n-1)}\frac{\tilde{\eta}_{0,i}}{n}\]
converge to $0$, $\Pb$-a.s. It follows that the process $\overline{X}_n:=\hat{Y}_{\overline{S}_n^{-1}}$ where
\[\overline{S}_n:=\sum_{x \in \Zb}\sum_{i=1}^{L^{\hat{Y}}(x,n-1)}\overline{\eta}_{x,i} \qquad \text{and} \qquad \overline{\eta}_{x,i}=\begin{cases} \eta_{x,i}, & \text{if } x < 0, \\ \tilde{\eta}_{\rho_x,i}, & \text{if } x\geq0, \end{cases} \]
obeys the same central limit theorem as $\hat{X}$. That is, we may replace the trap at the origin with the slightly different trap which corresponds to the branch at the root of the GW-tree, and still obtain a central limit theorem. 

Define the time spent above $0$ by $\overline{X}$ and the associated limiting Brownian motion $B$ as
\[A_t^{\overline{X}}:=\int_0^t\ind_{\{\overline{X}_s\geq 0\}}\d s \qquad \text{and} \qquad A_t^{B}:=\int_0^t\ind_{\{B_s\geq 0\}}\d s. \]
By substitution we then have that
\[\frac{A_{tn}^{\overline{X}}}{n}\; = \; \int_0^t\ind_{\{\overline{X}_{rn}/n^{1/2}\geq 0\}}\d r \; \rightarrow \; A_t^{B} \]
since $\lim_{\varepsilon \rightarrow 0^+}\int_0^t\ind_{\{B_s\in[-\varepsilon,\varepsilon]\}}\d s =0$. 
In particular, $(\overline{X}_{nt}n^{-1/2}, A_t^{\overline{X}}n^{-1})_{t\geq 0}$ converges to $(B_t, A_t^B)_{t\geq 0}$. 
Recall that $\tilde{X}$ is the projection of $X$ onto the backbone. By definition of $A_t^{\overline{X}}$, we have that $|\tilde{X}_{nt}|n^{-1/2}=\overline{X}_{(A_{tn}^{\overline{X}})^{-1}}n^{-1/2}$ which converges in distribution to $B_{(A_t^B)^{-1}}$ by \cite[Theorem 13.2.1 \& Corollary 13.6.4]{wh}. Furthermore, $B_{(A_t^B)^{-1}}$ is equal in distribution to $|B_t|$ hence we have the desired convergence result by Lemma \ref{proBack}.

Using Corollary \ref{SubTreeSpd} and taking the limit as $\beta \rightarrow 1^+$ we have that 
\[\frac{\nu_\beta}{\beta-1} \rightarrow \frac{\mu(1-\mu)}{2\sigma^2} =\frac{1}{2\Eb[\eta_0]}\]
by Lemma \ref{spd}, which completes the proof.
\end{proof}
\end{lem}
\subsection{An annealed functional central limit theorem}\label{SACLT}
We now prove an annealed functional central limit theorem for the biased walk on the subcritical GW-tree conditioned to survive. By Lemmas \ref{proBack} and \ref{Back} it will suffice to show the result holds for the corresponding randomly trapped random walk. We obtain the result by using the annealed invariance principle Theorem \ref{AIP}. That is, we show conditions on the tree and the bias which ensure that $\Eb[\eta_0^2]<\infty$. 

In order to show this we will use a decomposition which counts the number of visits to each vertex. The return probability given in Lemma \ref{returnProb} will be important. For $z_1,z_2,z_3 \in T$ write 
\[q_{z_1}(z_2,z_3):=\Pt^{T}_{z_1}(\tau^+_{z_2}<\tau^+_{z_3})\] 
to be the probability that the walk started from $z_1$ hits $z_2$ before $z_3$. We also require a similar relation for a walk on a tree. Let $T_{x,y}$ denote a tree with root $\rho$ in which every vertex has a single offspring except the vertices $w,x,y$ where $w$ has two offspring and $x,y$ have none. Denote these offspring $w_x,w_y$ then let $x,y$ be a descendants of $w_x,w_y$ respectively (possibly $w_x,w_y$). For vertices $z_1,z_2,z_3,z_4$ write 
\[q_{z_1}(z_2,\{z_3,z_4\}):=\Pt^{T_{x,y}}_{z_1}(\tau^+_{z_2}<\tau^+_{z_3}\land\tau^+_{z_4})\]
as the probability that the walk started at $z_1$ reaches $z_2$ before $z_3$ or $z_4$ by a $\beta$-biased walk on $T_{x,y}$.

\begin{figure}[H]
\centering
 \includegraphics[scale=0.9]{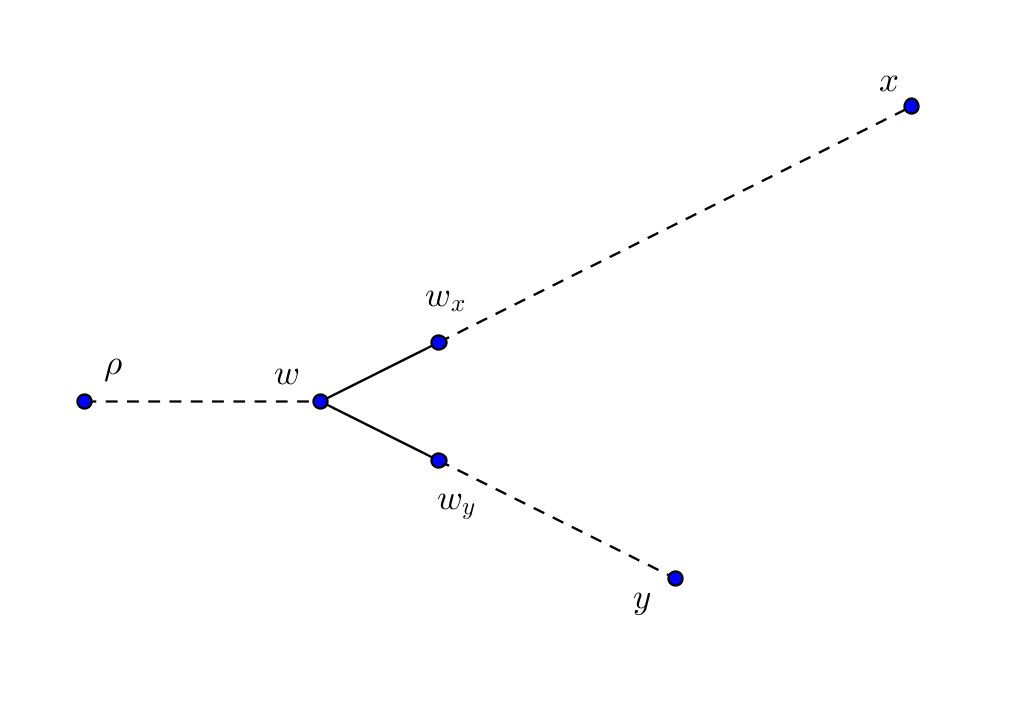} 
\caption{The tree $T_{x,y}$ with single branching point $w$ and extremal points $\rho,x,y$.}\label{Txy}
\end{figure}

Lemma \ref{retProb} gives the probability that the walk started at $w$ reaches $\rho$ before $x$ or $y$. Alternatively, this can be shown by comparing with an electrical network with conductances $\beta^k$ between vertices in generations $k,k+1$ and then using network reduction (see, for example, \cite[Chapter 2]{lype}).
\begin{lem}\label{retProb}
For any $T_{x,y}$,
\[ q_w(\rho,\{x,y\})  = \frac{(\beta^{|y|-|w|}-1)(\beta^{|x|-|w|}-1)}{2\beta^{|y|+|x|-|w|}-\beta^{|y|+|x|-2|w|}-\beta^{|x|}-\beta^{|y|}+1}.\]
 \begin{proof}
Write $w_\rho$ as the parent of $w$ then
\begin{flalign*}
 q_w(\rho,\{x,y\}) & = \frac{1}{2\beta+1}q_{w_\rho}(\rho,\{x,y\})+\frac{\beta}{2\beta+1}q_{w_x}(\rho,\{x,y\})+\frac{\beta}{2\beta+1}q_{w_y}(\rho,\{x,y\}) \\
 q_{w_\rho}(\rho,\{x,y\}) & =  q_w(\rho,\{x,y\})q_{w_\rho}(w,\rho)+q_{w_\rho}(\rho,w) \\
 q_{w_x}(\rho,\{x,y\}) & =  q_w(\rho,\{x,y\})q_{w_x}(w,x) \\
 q_{w_y}(\rho,\{x,y\}) & =  q_w(\rho,\{x,y\})q_{w_y}(w,y).
\end{flalign*}
Combining these gives us that
\begin{flalign*}
 q_w(\rho,\{x,y\}) & = \frac{q_w(\rho,\{x,y\})}{2\beta+1}\left(q_{w_\rho}(w,\rho)+\beta q_{w_x}(w,x)+\beta q_{w_y}(w,y)\right)+\frac{q_{w_\rho}(\rho,w)}{2\beta+1} \\
 & = \frac{q_{w_\rho}(\rho,w)}{2\beta+q_{w_\rho}(\rho,w)-q_{w_x}(w,x)-q_{w_y}(w,y)} \\
 & = \frac{\frac{\beta-1}{\beta^{|w|}-1}}{2\beta+\frac{\beta-1}{\beta^{|w|}-1}-\frac{\beta^{|x|-|w|}-\beta}{\beta^{|x|-|w|}-1}-\frac{\beta^{|y|-|w|}-\beta}{\beta^{|y|-|w|}-1}}
\end{flalign*}
by Lemma \ref{returnProb}. Rearranging gives the result.
 \end{proof}
\end{lem}

Let $T$ be a fixed tree and $(X_n)_{n\geq 1}$ a $\beta$-biased walk on $T$. For $x \in T$ let 
\[v_x:=\sum_{k=1}^{\tau_\rho^+}\ind_{\{X_k=x\}}\]
denote the number of visits to $x$ before returning to the root. Then $\tau_\rho^+=\sum_{x \in T}v_x$ and 
\begin{flalign}\label{ETrho2}
\Et_\rho^{T}\left[(\tau^+_\rho)^2\right]=\sum_{x,y \in T}\Et^{T}_\rho[v_xv_y].
\end{flalign}
For any $x,y \in T$ there exists a unique vertex $w_{x,y}$ which is the closest ancestor of both $x$ and $y$. We will often write $w$ instead of $w_{x,y}$ when it is clear to which vertices we are referring. Moreover \[\Et^{T}_\rho[v_xv_y]=\Pt^{T}_\rho(\tau^+_{w_{x,y}}<\tau^+_\rho)\Et^{T}_{w_{x,y}}[v_xv_y]\]
where, by comparison with a simple biased random walk on $\Zb$, we have that $\Pt^{T}(\tau^+_{w_{x,y}}<\tau^+_\rho)\in[ 1-\beta^{-1},1]$. We now prove a bound on $\Et^{T}_w[v_xv_y]$ following a similar method to that used in \cite{ke} for the unbiased case. Recall that $c(x)$ is the set of children of $x$ in $T$. 

\begin{lem}\label{vxvy}
For $\beta>1$, there exists a constant $C_\beta$ such that for any finite tree $T$,
\[\Et^{T}_\rho[v_xv_y]\leq\Et^{T}_w[v_xv_y]\leq C_\beta (|c(x)|\beta+1)(|c(y)|\beta+1)\beta^{|x|+|y|}. \]
\begin{proof}
When $w=\rho$ at least one of $x$ and $y$ is never reached therefore $v_xv_y=0$ and we may assume $|w|\geq 1$. There are now three cases to consider; these are:
\begin{enumerate}
 \item\label{xy} $x=y=w_{x,y}$;
 \item\label{xw} $x=w_{x,y}\neq y$;
 \item\label{ws} $x \neq w_{x,y}\neq y$.
\end{enumerate}
In case \ref{xy} we have that $v_x$ is geometrically distributed with termination probability $q_x(\rho,x)$ therefore 
\[\Et^{T}_{w_{x,y}}[v_xv_y]=\Et^{T}_{x}[v_x^2]=\frac{q_x(x,\rho)+1}{q_x(\rho,x)^2}.\]
For $x \notin c(\rho)$ we have that $\beta/(1+\beta)\leq q_x(x,\rho)  \leq 1$ and by Lemma \ref{returnProb} 
\[q_x(\rho,x)=\frac{1-\beta^{-1}}{(|c(x)|\beta+1)(\beta^{|x|-1}-\beta^{-1})}.\] 
We therefore have that 
\[\Et^{T}_{x}[v_x^2]\leq C_\beta (|c(x)|\beta+1)^2\beta^{2|x|}.\] 

In case \ref{xw}, the number of visits to $x$ from $x$ is geometrically distributed as in case \ref{xy}. For each visit to $x$ (except the last) the walk reaches $y$ before returning to $x$ with probability $q_x(y,x)/q_x(x,\rho)$ since, due to the tree structure, the walk cannot move from $\rho$ to $y$ without hitting $x$. From $y$, the walk returns to $y$ a geometric number of times before returning to $x$. More specifically,
\[\Et^{T}_{w_{x,y}}[v_xv_y]=\Et^{T}_{x}[v_xv_y]=\sum_{j=1}^\infty jq_x(\rho,x)q_x(x,\rho)^{j-1} \Et^{T}_x[v_y|v_x=j] \]
where, conditional on the event $\{v_x=j\}$, we have that $v_y$ is equal in distribution to the sum of
$B_{x,y}^j \sim Bin(j-1,q_x(y,x)/q_x(x,\rho))$ independent geometric random variables $G_{x,y}^i \sim Geo(q_y(x,y))$. Under $\Pt^{T}$ the number of excursions are independent therefore 
\[\Et^{T}_x[v_y|v_x=j]=(j-1)\frac{q_x(y,x)}{q_x(x,\rho)}\cdot\frac{1}{q_y(x,y)}.\]
We therefore have that 
\begin{flalign}\label{xwEx}
 \Et^{T}_{w_{x,y}}[v_xv_y] 
 & = \frac{q_x(y,x)q_x(\rho,x)}{q_x(x,\rho)q_y(x,y)}\sum_{j=1}^\infty j(j-1)q_x(x,\rho)^{j-1} \notag\\
 & =\frac{q_x(y,x)q_x(\rho,x)}{q_x(x,\rho)q_y(x,y)}\cdot \frac{2q_x(x,\rho)}{q_x(\rho,x)^3} \notag\\
  & =\frac{2q_x(y,x)}{q_y(x,y)q_x(\rho,x)^2}.
\end{flalign}
Using Lemma \ref{returnProb} we then have that
\begin{flalign*}
 q_x(y,x) & = \frac{\beta}{|c(x)|\beta+1}\cdot\frac{1-\beta^{-1}}{1-\beta^{|x|-|y|}}, \\
 q_y(x,y) & = \frac{1}{|c(y)|\beta+1}\cdot \frac{\beta-1}{\beta^{|y|-|x|}-1}, \\
 q_x(\rho,x) & = \frac{1}{|c(x)|\beta+1}\cdot \frac{\beta-1}{\beta^{|x|}-1}.
\end{flalign*}
Combining these with (\ref{xwEx}) we have that \[\Et^{T}_{w_{x,y}}[v_xv_y]\leq C_\beta (|c(x)|\beta+1)(|c(y)|\beta+1)\beta^{|x|+|y|}.\]

In case \ref{ws}, started from $w_{x,y}$, the walk reaches either $x$ or $y$ before returning to $\rho$ with probability $q_{w_{x,y}}(\{x,y\},\rho)$. From $x$ the walk has a geometric number of returns to $x$ before returning to $w_{x,y}$. Moreover, from $x$, the walk must return to $w_{x,y}$ before reaching either $\rho$ or $y$ by definition of $w_{x,y}$. The same also holds switching $x$ and $y$. Letting
\[\overline{q}_w(x,y)=\Pt^{T}_w\left(\tau^+_x<\tau^+_y|\tau^+_{\{x,y\}}<\tau^+_\rho\right) \qquad \text{and} \qquad \overline{q}_w(y,x)=\Pt^{T}_w\left(\tau^+_y<\tau^+_x|\tau^+_{\{x,y\}}<\tau^+_\rho\right) \]  
we then have that
\begin{flalign}\label{3vxvy}
 \Et^{T}_w[v_xv_y] & = \sum_{j=0}^\infty q_w(\{x,y\},\rho)^jq_w(\rho,\{x,y\}) \sum_{k=0}^j \overline{q}_w(x,y)^k\overline{q}_w(y,x)^{j-k}\binom{j}{k}\frac{k(j-k)}{q_x(w,x)q_y(w,y)}
\end{flalign}
since $q_x(w,x)^{-1}$ is the expected number of visits to $x$ (started from $x$) before returning to $w$ (and similarly for $y$) which are independent.
\begin{flalign*}
 & \sum_{k=0}^j \overline{q}_w(x,y)^k\overline{q}_w(y,x)^{j-k}\binom{j}{k}\frac{k(j-k)}{q_x(w,x)q_y(w,y)} \\
 & \qqqqquad = \frac{1}{q_x(w,x)q_y(w,y)}\sum_{k=1}^{j-1} \overline{q}_w(x,y)^k\overline{q}_w(y,x)^{j-k}\frac{j!}{(k-1)!(j-k-1)!} \\
 & \qqqqquad =  j(j-1)\frac{\overline{q}_w(x,y)\overline{q}_w(y,x)}{q_x(w,x)q_y(w,y)}\sum_{l=0}^{j-2} \overline{q}_w(x,y)^l\overline{q}_w(y,x)^{j-2-l}\frac{(j-2)!}{l!(j-2-l)!} \\
 & \qqqqquad =  j(j-1)\frac{\overline{q}_w(x,y)\overline{q}_w(y,x)}{q_x(w,x)q_y(w,y)}.
\end{flalign*}
Substituting back into (\ref{3vxvy}) it follows that
\begin{flalign*}
 \Et^{T}_w[v_xv_y] & = \frac{\overline{q}_w(x,y)\overline{q}_w(y,x)q_w(\rho,\{x,y\})}{q_x(w,x)q_y(w,y)}\sum_{j=0}^\infty j(j-1)q_w(\{x,y\},\rho)^j \\
 & = \frac{\overline{q}_w(x,y)\overline{q}_w(y,x)q_w(\rho,\{x,y\})}{q_x(w,x)q_y(w,y)} \cdot \frac{2q_w(\{x,y\},\rho)^2}{q_w(\rho,\{x,y\})^3} \\
 & = \frac{2q_w(\{x,y\},\rho)^2\overline{q}_w(x,y)\overline{q}_w(y,x)}{q_w(\rho,\{x,y\})^2q_x(w,x)q_y(w,y)}. 
\end{flalign*}
The terms in the numerator can all be bounded below by half of the escape probability $1-\beta^{-1}$ therefore we gain nothing using their exact expressions and bound them above by $1$. Using Lemmas \ref{returnProb} and \ref{retProb} for the other terms we have that
\begin{flalign*}
 q_w(\rho,\{x,y\}) & = \frac{(\beta^{|y|-|w|}-1)(\beta^{|x|-|w|}-1)}{2\beta^{|y|+|x|-|w|}-\beta^{|y|+|x|-2|w|}-\beta^{|x|}-\beta^{|y|}+1}, \\
 q_x(w,x) & = \frac{1}{|c(x)|\beta+1} \cdot \frac{\beta-1}{\beta^{|x|-|w|}-1}, \\
 q_y(w,y) & = \frac{1}{|c(y)|\beta+1} \cdot \frac{\beta-1}{\beta^{|y|-|w|}-1}.
\end{flalign*}
Since $|y|\geq 1$ we have that $\beta^{|y|}\geq 1$ therefore
\begin{flalign*}
q_w(\rho,\{x,y\}) \geq \frac{(\beta^{|y|-|w|}-1)(\beta^{|x|-|w|}-1)}{2\beta^{|y|+|x|-|w|}}
\end{flalign*}
and
\[\Et^{T}_w[v_xv_y]\leq\frac{2}{q_w(\rho,\{x,y\})^2q_x(w,x)q_y(w,y)}\leq C_\beta (|c(x)|\beta+1)(|c(y)|\beta+1)\beta^{|x|+|y|}. \]
\end{proof}
\end{lem}

Recall from (\ref{Nlab}) that $N$ is the number of return times to the root $\rho$ before reaching its unique ancestor $\oR$ and $\zeta_k$ is the duration of the $k\th$ such excursion. Letting $\oT=\oT_0$, by (\ref{eta1}) we have that
  \[\Et^{\omega}\left[\eta_0^2\right]=\Et^{\oT}\left[\sum_{k=1}^N\zeta_k^2+\sum_{k=1}^N\sum_{j\neq k}\zeta_k\zeta_j + 2 \sum_{k=1}^N\zeta_k+1\right].\]
We want to show that $\Er\left[  \Et^{\omega}\left[\eta_0^2\right]\right]<\infty$. Lemma \ref{AnBnd1} shows that this can be reduced to showing that the expected value of the first sum in the quenched expectation is finite.

\begin{lem}\label{AnBnd1}
    If $\beta^2\mu<1$ and $\Er[\xi^3]<\infty$ then \[\Er\left[\Et^{\oT}\left[\sum_{k=1}^N\sum_{j\neq k}\zeta_k\zeta_j + 2 \sum_{k=1}^N\zeta_k+1\right]\right]<\infty.\]
    \begin{proof}
       
By (\ref{eta1}) we have that \[\Et^{\oT}\left[2 \sum_{k=1}^N\zeta_k+1\right]< 2\Et^{\omega}[\eta_0]\]
which has finite expectation under $\Pr$ by Lemma \ref{spd} since $\beta^2\mu<1$ implies that $\beta\mu<1$.
  
The variable $N$ is geometrically distributed with termination probability $1-p_{ex}$; that is, 
\begin{flalign}\label{pex}
\Pt^{\oT}(N=k)=p_{ex}^k(1-p_{ex}) \quad \text{where} \quad p_{ex}=\frac{\beta(\xi^*-1)}{\beta\xi^*+1}
\end{flalign}
and $\xi^*+1$ is the number of neighbours attached to $\rho$ in $\oT$. We then have that 
\begin{flalign}
 \Et^{\oT}\left[\sum_{k=1}^N\sum_{j\neq k}\zeta_k\zeta_j\right] & = \sum_{n=0}^\infty \Pt^\oT(N=n)\sum_{k=1}^n\sum_{j\neq k} \Et^\oT[\zeta_k\zeta_j] \notag\\
 & =  \Et^\oT[\zeta_1]^2\sum_{n=0}^\infty n(n-1)p_{ex}^n(1-p_{ex}) \notag\\
 & =  2\left(\frac{\beta}{\beta+1}\right)^2\Et^\oT[\zeta_1]^2(\xi^*-1)^2 \label{pexsq}
\end{flalign}
by independence of excursion times under $\Pt^\oT$. 

Let $\oT^\circ$ denote the tree $\oT$ without the ancestor of the root $\oR$ and $\Tc$ be an $f$-GW-tree. Write $Z_n$ and $Z_n^\Tc$ to be the $n\th$ generation sizes of $\oT^\circ$ and $\Tc$ respectively then for $j=1,2,...$ let $Z_n^{\Tc,j}$ be independent copies of $Z_n^\Tc$. By the construction of $\oT^\circ$ using GW-trees we have that
\[ \sum_{l=1}^\infty l^2 \Pr\left(Z_k=l|Z_1=i\right)=\Er\left[(Z_{k-1}^{\Tc})^2|Z_0^{\Tc}=i\right]=\Er\left[\left(\sum_{j=1}^i Z_{k-1}^{\Tc,j}\right)^2\right]=i\Er\left[(Z_{k-1}^{\Tc})^2\right]+i(i-1)\Er\left[Z_{k-1}^{\Tc}\right]^2.\]

Using this with Lemma \ref{gensize}, for $j\geq k\geq 1$ we have that
\begin{flalign}
 \Er\left[Z_kZ_j\right] 
 & = \sum_{i=1}^\infty\Pr\left(\xi^*-1=i\right) \sum_{l=1}^\infty l \Pr\left(Z_k=l|Z_1=i\right)\Er\left[Z_j|Z_k=l\right] \notag\\
 & = \sum_{i=1}^\infty\Pr\left(\xi^*-1=i\right) \sum_{l=1}^\infty l^2 \Pr\left(Z_k=l|Z_1=i\right)\Er\left[Z_{j-k}^{\Tc}\right] \notag\\
 & = \mu^{j-k} \sum_{i=1}^\infty\Pr\left(\xi^*-1=i\right) \left(i\Er\left[(Z_{k-1}^{\Tc})^2\right]+i(i-1)\Er\left[Z_{k-1}^{\Tc}\right]^2\right) \notag\\
 & \leq C\mu^j\Er[(\xi^*-1)^2]. \label{ZkZj}
\end{flalign}

Using (\ref{pexsq}) and the formula (\ref{return}) for the expected time spent in a tree we have that 
\begin{flalign*}
 \Er\left[\Et^{\oT}\left[\sum_{k=1}^N\sum_{j\neq k}\zeta_k\zeta_j\right]\right] & = C_\beta \Er\left[(\xi^*-1)^2\Et^\Tc[\zeta_1]^2\right] \\
 & = C_\beta \Er\left[\left(\sum_{k\geq 1} \beta^kZ_k\right)^2\right] \\
 & = C_\beta\sum_{k\geq 1}\beta^{2k}\Er[Z_k^2] + 2C_\beta\sum_{k\geq 1}\sum_{j>k}\beta^{k+j}\Er[Z_kZ_j]. 
 \end{flalign*}
By (\ref{ZkZj}) we then have that
 \[\sum_{k\geq 1}\sum_{j>k}\beta^{k+j}\Er[Z_kZ_j]\;\leq\; C\Er[(\xi^*-1)^2]\sum_{k\geq 1}\beta^k\sum_{j>k}(\beta\mu)^j\;\leq\; C\Er[(\xi^*-1)^2]\sum_{k\geq 1}(\beta^2\mu)^k\]
 and
 \[\sum_{k\geq 1}\beta^{2k}\Er[Z_k^2] \; \leq \; C\Er[(\xi^*-1)^2]\sum_{k\geq 1}(\beta^2\mu)^k. \]
Each of these terms is finite since $\beta^2\mu<1$ and $\Er[\xi^3]<\infty$ where we recall from the definition of the size-biased distribution that $\Er[(\xi^*-1)^2]\leq C\Er[\xi^3]$.
    \end{proof}
  \end{lem}

In order to show that $\Eb[\eta_0^2]<\infty$ it remains to prove Lemma \ref{AnBnd2} which follows similarly to Lemma \ref{AnBnd1} with the use of Lemma \ref{vxvy}.

\begin{lem}\label{AnBnd2}
 If $\beta^2\mu<1$ and $\Er[\xi^3]<\infty$ then 
 \[\Er\left[\Et^{\oT}\left[\sum_{k=1}^N\zeta_k^2\right]\right]<\infty.\]
 \begin{proof}
Recall that $\oT^\circ$ denotes the tree $\oT$ without the ancestor of the root $\oR$. Since the separate excursions are independent under $\Pt^{\oT}$ and $N$ is geometrically distributed we have that 
\[\Er\left[\Et^{\oT}\left[\sum_{k=1}^N\zeta_k^2\right]\right]=\Er\left[\left(\frac{\beta(\xi^*-1)}{\beta+1}\right)\Et^{\oT^\circ}_\rho\left[\zeta_1^2\right]\right].\]
Labelling $\rho_{1},...,\rho_{\xi^*-1}$ as the neighbours of $\rho$ in $\oT^\circ$, and $\oT_{\!\!\rho_k}$ as the tree consisting of $\rho, \rho_k$ and the descendants of $\rho_k$ we have that 
\begin{flalign*}
 \Et^{\oT^\circ}_\rho\left[\zeta_1^2\right] & =\sum_{k=1}^{\xi^*-1} \frac{\Et^{\oT_{\!\!\rho_k}}_\rho\left[(\tau_\rho^+)^2\right]}{\xi^*-1}
\end{flalign*}
when $\xi^*\neq 1$ and $0$ otherwise. Moreover, it then follows that
\begin{flalign*}
 \Er\left[\Et^{\oT}\left[\sum_{k=1}^N\zeta_k^2\right]\right] & \leq 2\Er\left[\sum_{k=1}^{\xi^*-1}\Et^{\oT_{\!\!\rho_{k}}}\left[(\tau^+_\rho)^2\right]\right] = 2\Er[\xi^*-1]\Er\left[\Et^{\oT_{\!\!\rho_{1}}}\left[(\tau^+_\rho)^2\right]\right]
\end{flalign*}
since the subtraps are independent. Since $\Er[\xi^*-1]\leq C\Er[\xi^2]<\infty$, it suffices to show that 
\[\Er\left[\Et^{\ooT}_{\rho}\left[(\tau^+_\rho)^2\right]\right]<\infty\]
where $\ooT$ is a tree (equal in distribution to $\oT_{\!\rho_1}$) with root $\rho$, single first generation vertex $\overrightarrow{\rho}$ and, under $\Pr$, the subtree rooted at $\overrightarrow{\rho}$ is a subcritical GW-tree with the original offspring distribution.

Recall that $\ooT_z$ denotes the descendent tree of $\ooT$ at $z$. By (\ref{ETrho2}) and Lemma \ref{vxvy} we have that
\begin{flalign*}
\Er\left[\Et^{\ooT}_{\rho}\left[(\tau^+_\rho)^2\right]\right]
\; = \;\Er\left[\sum_{x,y \in \ooT}\Et^{\ooT}_\rho[v_xv_y]\right] 
\;\leq \;C_\beta\Er\left[\left(\sum_{x\in \ooT}\left(|c(x)|\beta+1\right)\beta^{|x|}\right)\left(\sum_{y\in \ooT}\left(|c(y)|\beta+1\right)\beta^{|y|}\right)\right].
\end{flalign*}

By collecting terms in the $k\th$ generation we have that 
\[\sum_{x\in \ooT}\left(|c(x)|\beta+1\right)\beta^{|x|}\;=\;1+\sum_{k\geq 1}Z_k^{\ooT}(\beta^k+\beta^{k-1}) \;\leq\; (1-\beta^{-1})\sum_{k\geq 0}Z_k^{\ooT}\beta^k\]
where $Z_k^{\ooT}$ is the size of the $k\th$ generation of $\ooT$. For $k\geq 0$ the tree $\ooT$ satisfies $Z^{\ooT}_{k+1}=Z_k$ for a GW-process $Z_k$ with $Z_0=1$; therefore, using that $\beta^2\mu<1$ and Lemma \ref{gensize}, we have that $\Er\left[Z_k^{\ooT}Z_j^{\ooT}\right] \leq C\mu^j$, for $j\geq k$. In particular, 
\begin{flalign*}
 \Er\left[\Et^{\ooT}_{\rho}\left[(\tau^+_\rho)^2\right]\right] 
  \;\leq\; C_\beta\sum_{k\geq 0}\beta^k\sum_{j\geq k}\beta^j\Er\left[Z_k^{\ooT}Z_j^{\ooT}\right]  
  \;\leq\; C_\beta\sum_{k\geq 0}\beta^k\sum_{j\geq k}(\mu\beta)^j 
  \;\leq\; C_{\beta,\mu}\sum_{k\geq 0}(\mu\beta^2)^k
  \;<\;\infty.
\end{flalign*}
\end{proof}
\end{lem}

By Lemmas \ref{AnBnd1} and \ref{AnBnd2} we have that $\Eb[\eta_0^2]<\infty$ therefore, by Lemmas \ref{proBack} and \ref{Back} and Theorem \ref{AIP}, we have the following annealed functional central limit theorem: 
\begin{cly}\label{AnnSecMom}
 If $\beta^2\mu<1$ and $\Er[\xi^3]<\infty$ then $\Eb[\eta_0^2]<\infty$ and, in particular, there exists $\varsigma^2<\infty$ such that
 \[B_t^n=\frac{|X_{nt}| -nt\nu_\beta}{\varsigma \sqrt{n}}\] converges in $\Pb$-distribution on $D(\Rb^+,\Rb)$ endowed with the Skorohod metric to a standard Brownian motion.  
\end{cly}

Recall that the expression (\ref{varsig}) for $\varsigma^2$ was given in Theorem \ref{AIP} in terms of the moments of the distance and time between regenerations. We can therefore use this to write the corresponding form in the GW-tree model as 
\[\varsigma^2=\frac{\Eb\left[\left(Y_{\kappa_2}-Y_{\kappa_1}-\nu_\beta\sum_{j=\kappa_1}^{\kappa_2-1}\eta_j\right)^2\right]}{\Eb[\eta_0]\Eb[\kappa_2-\kappa_1]}\]
where $\kappa_j$ are the regeneration times of the walk $Y$.
We now show that both of the conditions $\beta^2\mu<1$ and $\Er[\xi^3]<\infty$ are necessary in order to apply Theorem \ref{AIP} and thus obtain an annealed functional central limit theorem for the walk on the subcritical GW-tree conditioned to survive. 
\begin{lem}\label{Nec}
 If $\beta^2\mu\geq 1$ or $\Er[\xi^3]=\infty$ then \[\Eb[\eta_1^2]\geq \Er\left[\Et^\oT[\eta_1]^2\right]=\infty.\]
 \begin{proof}
 Recall that $\eta_1$ is the first hitting time of $\oR$ by $W_n$ started from the root $\rho$ in $\oT$. With positive probability $\rho$ has neighbours other than $\oR$ and the walk moves to one on the first step. Until returning to $\rho$ the walk is equal in distribution to a $\beta$-biased random walk on an $f$-GW-tree conditioned to have a single first generation vertex. In particular, it suffices to show that for a $\beta$-biased walk 
 \[\Er\left[\Et^{\Tc}_\rho\left[\tau^+_\rho\right]^2|Z_1=1\right]=\infty\]
 where $\Tc$ is an $f$-GW-tree rooted at $\rho$. Using the formula for the expected time spent in a tree (\ref{return}) we have that 
  \begin{flalign*}
   \Er\left[\Et^{\Tc}_\rho\left[\tau^+_\rho\right]^2|Z_1=1\right] 
   = \frac{4}{\beta^2}\Er\left[\left(\sum_{k\geq 1}\beta^kZ_k\right)^2|Z_1=1\right]  
    \geq \frac{4}{\beta^2}\sum_{k\geq 1}\beta^{2k}\Er[Z_k^2|Z_1=1].
  \end{flalign*}
Since $Z_k$ takes nonnegative values in $\Zb$ we have that 
\[\Er\left[Z_k^2|Z_1=1\right]\geq \Er\left[Z_k|Z_1=1\right] =\mu^{k-1}\]
by statement \ref{EZn} of Lemma \ref{gensize}. We therefore have that 
\[\Er\left[\Et^{\Tc}_\rho[\tau^+_\rho]^2|Z_1=1\right]\geq c\sum_{k\geq 1}(\beta^2\mu)^k\]
which is infinite if $\beta^2\mu\geq 1$. 

The first hitting time of $\oR$ is at least the number of visits to the offspring of $\rho$. From $\rho$, the walk takes a geometric number of visits (with termination probability $1-p_{ex}$, see (\ref{pex})) to these vertices before reaching $\oR$. Using properties of geometric random variables we then have that
\begin{flalign*}
 \Er\left[\Et^\oT_\rho[\eta_1]^2\right] \;\geq\; \Er\left[\left(\frac{(\xi^*-1)\beta}{\beta+1}\right)^2\right] \;  \geq \; c\left(\Er[(\xi^*)^2]-1\right) \;  =\; c \left(\mu^{-1}\Er[\xi^3]-1\right).
\end{flalign*}
 \end{proof}
\end{lem}

\subsection{A quenched central limit theorem}\label{SQCLT}
We now prove a quenched central limit theorem for the biased walk on the subcritical GW-tree conditioned to survive. As in the annealed case, by Lemmas \ref{proBack} and \ref{Back} it will suffice to show the result holds for the corresponding randomly trapped random walk and we obtain the result by using Theorem \ref{QCLT}. As in the randomly trapped random walk case, he value of $\vartheta$ is known to be $\varsigma \nu_\beta^{3/2}$ where $\nu_\beta$ and $\varsigma$ are given in Corollary \ref{SubTreeSpd} and \ref{AnnSecMom} respectively. Notice that, under the assumptions of Proposition \ref{enqusq}, $\nu_\beta t-\Gc^{\Tc^*}(t)$ is a sum of i.i.d.\ centred random variables with positive, finite variance. It therefore follows that this expression converges in distribution with respect to $\Pr$ to a Gaussian random variable. In particular, this means that the environment dependent centring is necessary.

\begin{prp}\label{enqusq}
 If $\beta^2\mu<1$ and $\Er[\xi^{3+\delta}]<\infty$ for some $\delta>0$ then 
 \[\Er\left[\Et^\oT\left[\eta_1\right]^{2+\varepsilon}\right]<\infty\]
 for some $\varepsilon>0$. Moreover, there exists $\vartheta>0$ such that for $\Pr$-a.e.\ $\Tc^*$ we have that \[\Pt^{\Tc^*}\left(\frac{X_t-\Gc^{\Tc^*}(t)}{\vartheta\sqrt{t}} \leq x\right) \rightarrow \Phi(x)\]
 uniformly in $x$ as $\nin$ where $\Phi(x)$ is given in Lemma \ref{HitCLT} and for $\nu_\beta$ as in Corollary \ref{SubTreeSpd} 
 \[\Gc^{\Tc^*}(t)=\nu_\beta t - \nu_\beta\sum_{k=1}^{\lfloor \nu_\beta t\rfloor}\frac{\beta+1}{\beta-1}\left(\Et^{\Tc^*}[\tilde{\eta}_{\rho_k,0}]-\Eb[\eta_0]\right).\]
 
 \begin{proof}
 By Proposition \ref{AnnSecMom} we have that $\Eb[\eta_0^2]<\infty$ therefore by the previous remark it suffices to show that for some $\varepsilon\in(0,1)$
 \[\Er\left[\Et^\oT\left[\eta_1\right]^{2+\varepsilon}\right]<\infty.\]
  
 Recall that \[N=\sum_{n=1}^{\tau^+_{\overleftarrow{\rho}}}\ind_{\{W_n=\rho\}}\]
 is the number of hitting times of the root $\rho$ before reaching $\overleftarrow{\rho}$ (for the walk started at $\rho$) and $\oT^\circ$ is the tree $\oT$ with $\overleftarrow{\rho}$ removed. Then 
 \[\Et^\oT_\rho\left[\eta_1\right]\;=\;1+N \Et^{\oT^\circ}_\rho\left[\tau^+_\rho\right]\;=\; 1+N \sum_{n\geq 1}\frac{Z_n}{Z_1}\beta^{n-1} \;\leq\; (N+1) \sum_{n\geq 1}\frac{Z_n}{Z_1}\beta^{n-1} \]
 by (\ref{return}) where $Z_n$ is the $n\th$ generation size of $\oT^\circ$ since the walk on $\oT^\circ$ is $\beta$-biased.
 
 For a fixed tree, $N$ is geometrically distributed with excursion probability $p_{ex}$ (see (\ref{pex})). By conditioning on $Z_1$ we therefore have that
 \begin{flalign*}
 \Er\left[\Et^\oT\left[\eta_1\right]^{2+\varepsilon}\right] & \leq \Er\left[\Er\left[(N+1)^{2+\varepsilon}|Z_1\right]\Er\left[\left(\sum_{n\geq 1}\frac{Z_n}{Z_1}\beta^{n-1}\right)^{2+\varepsilon}\Big|Z_1\right]\right] \\
 & \leq C\Er\left[Z_1^{2+\varepsilon}\Er\left[\left(\sum_{n\geq 1}\frac{Z_n}{Z_1}\beta^{n-1}\right)^{2+\varepsilon}\Big|Z_1\right]\right] \\
 & = C\Er\left[\left(\sum_{n\geq 1}Z_n\beta^{n-1}\right)^{2+\varepsilon}\right]. 
 \end{flalign*}
 
 We can write \[Z_n=\sum_{j=1}^{Z_1}Z_{n-1}^{(j)}\]
 where $Z^{(j)}$ are independent GW-processes therefore by convexity
 \begin{flalign*}
  \Er\left[\left(\sum_{n\geq 1}Z_n\beta^{n-1}\right)^{2+\varepsilon}\right] & = \Er\left[Z_1^{2+\varepsilon}\left(\sum_{j=1}^{Z_1}\sum_{n\geq 1}\frac{Z_{n-1}^{(j)}}{Z_1}\beta^{n-1}\right)^{2+\varepsilon}\right] \\
  & \leq \Er\left[Z_1^{1+\varepsilon}\sum_{j=1}^{Z_1}\left(\sum_{n\geq 1}Z_{n-1}^{(j)}\beta^{n-1}\right)^{2+\varepsilon}\right] \\
  & = \Er[(\xi^*-1)^{2+\varepsilon}]\Er\left[\left(\sum_{n\geq 1}Z_{n-1}^{(1)}\beta^{n-1}\right)^{2+\varepsilon}\right].
 \end{flalign*}
By the assumptions of the theorem we have that $\Er[(\xi^*-1)^{2+\varepsilon}]\leq \mu^{-1}\Er[\xi^{3+\varepsilon}]<\infty$ whenever $\varepsilon<\delta$ thus it suffices to show that 
\[\Er\left[\left(\sum_{n\geq 0}Z_{n}\beta^{n}\right)^{2+\varepsilon}\right]<\infty\]
where $Z_n$ now denotes the $n\th$ generation size of an $f$-GW-process.
 
 For $\varepsilon<\delta$, by conditioning on the height $\Hc:=\max\{n\geq 0: Z_n>0\}$ of the tree we have that
\begin{flalign}\label{conH}
 \Er\left[\left(\sum_{n\geq 0}Z_n\beta^n\right)^{2+\varepsilon}\right]  &  \leq \Er\left[\beta^{(2+\varepsilon)\Hc}\Er\left[\left(\sum_{n\geq 0}Z_n\right)^{2+\varepsilon}\Big|\Hc\right]\right] \notag\\
 & \leq \Er\left[\beta^{(2+\varepsilon)\Hc}(\Hc+1)^{2+\varepsilon}\Er\left[\max_{n\leq \Hc}Z_n^{2+\varepsilon}\Big|\Hc\right]\right]\notag \\
 & \leq \Er\left[\beta^{(2+\varepsilon)\Hc}(\Hc+1)^{2+\varepsilon}\sum_{n=0}^{\Hc}\Er\left[Z_n^{2+\varepsilon}\Big|\Hc\right]\right] \notag\\
 & = \sum_{n=0}^\infty \Er\left[\beta^{(2+\varepsilon)\Hc}(\Hc+1)^{2+\varepsilon}Z_n^{2+\varepsilon}\right] \notag\\
 & = \sum_{n=0}^\infty\sum_{j=1}^\infty j^{2+\varepsilon}\Pr(Z_n=j)\Er\left[\beta^{(2+\varepsilon)\Hc}(\Hc+1)^{2+\varepsilon}|Z_n=j\right]. 
\end{flalign}
From (\ref{cmu}) we have that $\Pr(\Hc\geq n) \sim c\mu^n$ therefore $\Pr(\Hc\geq n)\leq C\mu^n$ for some constant $C$ hence
\begin{flalign*}
 \Er\left[\beta^{(2+\varepsilon)\Hc}(\Hc+1)^{2+\varepsilon}|Z_n=j\right] & = \Er\left[\beta^{(2+\varepsilon)(\Hc+n+1)}(\Hc+n)^{2+\varepsilon}|Z_0=j\right] \\
 & = \sum_{i=1}^\infty \beta^{(2+\varepsilon)(i+n)}(i+n+1)^{2+\varepsilon}\Pr(\Hc=i|Z_0=j) \\
 & \leq \sum_{i=1}^\infty \beta^{(2+\varepsilon)(i+n)}(i+n+1)^{2+\varepsilon}\Pr(\Hc\geq i|Z_0=j) \\
 & \leq C\sum_{i=1}^\infty \beta^{(2+\varepsilon)(i+n)}(i+n+1)^{2+\varepsilon}j\mu^i \\
 & \leq Cj\beta^{(2+\varepsilon)n}(n+2)^{2+\varepsilon}\sum_{i=1}^\infty i^{2+\varepsilon}(\beta^{2+\varepsilon}\mu)^i.
\end{flalign*}
Since $\beta^2\mu<1$ we can choose $\varepsilon>0$ sufficiently small so that $\beta^{2+\varepsilon}\mu<1$ therefore
\[\sum_{i=1}^\infty i^{2+\varepsilon}(\beta^{2+\varepsilon}\mu)^i<\infty.\]

Substituting back into (\ref{conH}) we have that
\begin{flalign}\label{EET2ep}
 \Er\left[\left(\sum_{n\geq 0}Z_n\beta^n\right)^{2+\varepsilon}\right]  & \leq C\sum_{n=0}^\infty\beta^{n(2+\varepsilon)}(n+2)^{2+\varepsilon}\sum_{j=1}^\infty j^{3+\varepsilon}\Pr(Z_n=j) \notag \\
  & = C\sum_{n=0}^\infty\beta^{n(2+\varepsilon)}(n+2)^{2+\varepsilon}\Er[Z_n^{3+\varepsilon}].
\end{flalign}

Using a telescoping sum we can write
\[Z_n=\mu^n+\sum_{k=0}^{n-1}(Z_{n-k}-\mu Z_{n-(k+1)})\mu^k,\]
therefore using convexity we have that
\begin{flalign}\label{Zn3ep}
 \Er[Z_n^{3+\varepsilon}] & = (n+1)^{3+\varepsilon}\Er\left[\left(\frac{\mu^n}{n+1}+\sum_{k=0}^{n-1}\frac{(Z_{n-k}-\mu Z_{n-(k+1)})\mu^k}{n+1}\right)^{3+\varepsilon}\right] \notag \\
 & \leq (n+1)^{3+\varepsilon}\Er\left[\frac{\mu^{n(3+\varepsilon)}}{n+1}+\sum_{k=0}^{n-1}\frac{\left((Z_{n-k}-\mu Z_{n-(k+1)})\mu^k\right)^{3+\varepsilon}}{n+1}\right] \notag \\
 & = (n+1)^{2+\varepsilon}\mu^{n(3+\varepsilon)}+(n+1)^{2+\varepsilon}\sum_{k=0}^{n-1}\mu^{k(3+\varepsilon)}\Er\left[\left(Z_{n-k}-\mu Z_{n-(k+1)}\right)^{3+\varepsilon}\right].
\end{flalign}

Let $\xi_j$ be independent copies of $\xi$ then using the Marcinkiewicz-Zygmund inequality and convexity we have that
\begin{flalign*}
 \Er\left[\left(Z_{n-k}-\mu Z_{n-(k+1)}\right)^{3+\varepsilon}\right] & = \Er\left[\Er\left[\left(\sum_{j=1}^{Z_{n-(k+1)}}(\xi_j-\mu)\right)^{3+\varepsilon}\Big| Z_{n-(k+1)}\right]\right] \\
 & \leq C\Er\left[\Er\left[\left(\sum_{j=1}^{Z_{n-(k+1)}}(\xi_j-\mu)^2\right)^{\frac{3+\varepsilon}{2}}\Big| Z_{n-(k+1)}\right]\right] \\
 & = C\Er\left[\Er\left[\left(\sum_{j=1}^{Z_{n-(k+1)}}\frac{(\xi_j-\mu)^2}{Z_{n-(k+1)}}\right)^{\frac{3+\varepsilon}{2}}\Big| Z_{n-(k+1)}\right]Z_{n-(k+1)}^{\frac{3+\varepsilon}{2}}\right] \\
 & \leq C\Er\left[\Er\left[\sum_{j=1}^{Z_{n-(k+1)}}\frac{|\xi_j-\mu|^{3+\varepsilon}}{Z_{n-(k+1)}}\Big| Z_{n-(k+1)}\right]Z_{n-(k+1)}^{\frac{3+\varepsilon}{2}}\right] \\
 & = C\Er\left[|\xi-\mu|^{3+\varepsilon}\right]\Er\left[Z_{n-(k+1)}^{\frac{3+\varepsilon}{2}}\right] \\
 & \leq C\Er\left[|\xi-\mu|^{3+\varepsilon}\right]\Er\left[Z_{n-(k+1)}^2\right]. 
\end{flalign*}

By Lemma \ref{gensize} we have that $\Er\left[Z_{n-(k+1)}^2\right]\leq C\mu^{n-(k+1)}$ where $C$ is independent of $n,k$ therefore substituting into (\ref{Zn3ep}) we have that 
\[\Er[Z_n^{3+\varepsilon}] \;\leq\; (n+1)^{2+\varepsilon}\mu^{n(3+\varepsilon)}+C(n+1)^{2+\varepsilon}\mu^n\sum_{k=0}^{n-1}\mu^{k(2+\varepsilon)} \;\leq\; C(n+1)^{2+\varepsilon}\mu^n.\]
Combining with (\ref{EET2ep}) we then have that \[\Er\left[\left(\sum_{n\geq 0}Z_n\beta^n\right)^{2+\varepsilon}\right]\leq C\sum_{n=1}^\infty(n+2)^{4+\varepsilon}(\beta^{2+\varepsilon}\mu)^n \] 
which is finite since we have chosen $\varepsilon>0$ sufficiently small so that $\beta^{2+\varepsilon}\mu<1$.
\end{proof}

\end{prp}

\section*{Acknowledgements}
I would like to thank my supervisor David Croydon for suggesting the problem, his support and many useful discussions. This work is supported by EPSRC as part of the MASDOC DTC at the University of Warwick. Grant No.\ EP/H023364/1.


\begin{thebibliography}{10}

\bibitem{ai}
E.~A{\"{\i}}d{\'e}kon.
\newblock Speed of the biased random walk on a {G}alton-{W}atson tree.
\newblock {\em Probab. Theory Related Fields}, 159(3-4):597--617, 2014.

\bibitem{arcacero}
G.~Ben~Arous, M.~Cabezas, J.~{\v{C}}ern{\'y}, and R.~Royfman.
\newblock Randomly trapped random walks.
\newblock {\em Ann. Probab.}, 43(5):2405--2457, 2015.

\bibitem{arfr}
G.~Ben~Arous and A.~Fribergh.
\newblock Biased random walks on random graphs.
\newblock {\em Probab. Statist. Phys. St. Petersburg}, 91:99, 2016.

\bibitem{arfrgaha}
G.~Ben~Arous, A.~Fribergh, N.~Gantert, and A.~Hammond.
\newblock Biased random walks on {G}alton-{W}atson trees with leaves.
\newblock {\em Ann. Probab.}, 40(1):280--338, 2012.

\bibitem{bosz}
E.~Bolthausen and A.~Sznitman.
\newblock On the static and dynamic points of view for certain random walks in
  random environment.
\newblock {\em Methods and Applications of Analysis}, 2002.

\bibitem{bouc}
J.~Bouchaud.
\newblock {Weak ergodicity breaking and aging in disordered systems}.
\newblock {\em {J. Phys. I}}, 2(9):1705--1713, 1992.

\bibitem{bopre}
A.~Bowditch.
\newblock A quenched central limit theorem for biased random walks on
  supercritical {G}alton-{W}atson trees.
\newblock (in preparation).

\bibitem{bowd}
A.~Bowditch.
\newblock Escape regimes of biased random walks on {G}alton-{W}atson trees.
\newblock {\em arXiv preprint arXiv:1605.05050}, 2016.

\bibitem{cewa}
J.~{\v{C}}ern{\'y} and T.~Wassmer.
\newblock Randomly trapped random walks on {$\mathbb{Z}^d$}.
\newblock {\em Stochastic Process. Appl.}, 125(3):1032--1057, 2015.

\bibitem{crfrku}
D.~A. Croydon, A.~Fribergh, and T.~Kumagai.
\newblock Biased random walk on critical {G}alton-{W}atson trees conditioned to
  survive.
\newblock {\em Probab. Theory Related Fields}, 157(1-2):453--507, 2013.

\bibitem{depeze}
A.~Dembo, Y.~Peres, and O.~Zeitouni.
\newblock Tail estimates for one-dimensional random walk in random environment.
\newblock {\em Comm. Math. Phys.}, 181(3):667--683, 1996.

\bibitem{dh}
D.~Dhar.
\newblock Diffusion and drift on percolation networks in an external field.
\newblock {\em J. of Phys.}, 17(5):L257, 1984.

\bibitem{du}
R.~Durrett.
\newblock {\em Probability: theory and examples}.
\newblock Cambridge Series in Statistical and Probabilistic Mathematics.
  Cambridge University Press, Cambridge, fourth edition, 2010.

\bibitem{foisne}
L.~Fontes, M.~Isopi, and C.~Newman.
\newblock Random walks with strongly inhomogeneous rates and singular
  diffusions: convergence, localization and aging in one dimension.
\newblock {\em Ann. Probab.}, 30(2):579--604, 2002.

\bibitem{frha}
A.~Fribergh and A.~Hammond.
\newblock Phase transition for the speed of the biased random walk on the
  supercritical percolation cluster.
\newblock {\em Comm. Pure Appl. Math.}, 67(2):173--245, 2014.

\bibitem{go}
I.~Goldsheid.
\newblock Simple transient random walks in one-dimensional random environment:
  the central limit theorem.
\newblock {\em Probab. Theory Related Fields}, 139(1-2):41--64, 2007.

\bibitem{heseve}
C.~R. Heathcote, E.~Seneta, and D.~Vere-Jones.
\newblock A refinement of two theorems in the theory of branching processes.
\newblock {\em Teor. Verojatnost. i Primenen.}, 12:341--346, 1967.

\bibitem{ja}
S.~Janson.
\newblock Simply generated trees, conditioned {G}alton-{W}atson trees, random
  allocations and condensation.
\newblock {\em Probab. Surv.}, 9:103--252, 2012.

\bibitem{ke}
H.~Kesten.
\newblock Subdiffusive behavior of random walk on a random cluster.
\newblock {\em Ann. Inst. H. Poincar\'e Probab. Statist.}, 22(4):425--487,
  1986.

\bibitem{kekosp}
H.~Kesten, M.~V. Kozlov, and F.~Spitzer.
\newblock A limit law for random walk in a random environment.
\newblock {\em Compos. Math.}, 30:145--168, 1975.

\bibitem{lypepeot}
R.~Lyons, R.~Pemantle, and Y.~Peres.
\newblock Conceptual proofs of {$L\log L$} criteria for mean behavior of
  branching processes.
\newblock {\em Ann. Probab.}, 23(3):1125--1138, 1995.

\bibitem{lypepe}
R.~Lyons, R.~Pemantle, and Y.~Peres.
\newblock Biased random walks on {G}alton-{W}atson trees.
\newblock {\em Probab. Theory Related Fields}, 106(2):249--264, 1996.

\bibitem{lype}
R.~Lyons and Y.~Peres.
\newblock {\em Probability on Trees and Networks}.
\newblock Cambridge University Press, New York, 2016.

\bibitem{mazy}
J.~Marcinkiewicz and A.~Zygmund.
\newblock Quelques th{\'e}oremes sur les fonctions ind{\'e}pendantes’.
\newblock {\em Fund. Math}, 29:60--90, 1937.

\bibitem{peze}
Y.~Peres and O.~Zeitouni.
\newblock A central limit theorem for biased random walks on {G}alton-{W}atson
  trees.
\newblock {\em Probab. Theory Related Fields}, 2008.

\bibitem{pe}
V.~V. Petrov.
\newblock {\em Sums of independent random variables}.
\newblock Springer-Verlag, New York-Heidelberg, 1975.

\bibitem{sz}
A.~Sznitman.
\newblock Slowdown estimates and central limit theorem for random walks in
  random environment.
\newblock {\em J. Eur. Math. Soc. (JEMS)}, 2(2):93--143, 2000.

\bibitem{szit}
A.~Sznitman.
\newblock On the anisotropic walk on the supercritical percolation cluster.
\newblock {\em Comm. Math. Phys.}, 240(1-2):123--148, 2003.

\bibitem{szze}
A.~Sznitman and M.~Zerner.
\newblock A law of large numbers for random walks in random environment.
\newblock {\em Ann. Probab.}, 27(4):1851--1869, 1999.

\bibitem{wh}
W.~Whitt.
\newblock {\em Stochastic-process limits}.
\newblock Springer-Verlag, New York, 2002.

\bibitem{ya}
A.~M. Yaglom.
\newblock Certain limit theorems of the theory of branching random processes.
\newblock {\em Doklady Akad. Nauk SSSR (N.S.)}, 56:795--798, 1947.

\bibitem{zi}
O.~Zindy.
\newblock Scaling limit and aging for directed trap models.
\newblock {\em Markov Process. Related Fields}, 15(1):31--50, 2009.

\end{thebibliography}
\end{document}